\documentclass[3p,review,fleqn,number,sort&compress,times]{elsarticle}

\usepackage{amssymb}
\usepackage{amsmath}
\usepackage{booktabs}
\usepackage{calc}
\usepackage{color}
\usepackage{enumitem}
\usepackage{esvect}
\usepackage{icomma}
\usepackage{pgfplots}
\usepackage{pifont}
\usepackage{siunitx}
\usepackage{stackengine}
\usepackage{subcaption}
\usepackage{tikz}
\usepackage{xfrac}
\usepackage{xspace}

\usepackage{accents}
\usepackage[colorlinks]{hyperref}
\usepackage[capitalize]{cleveref}

\usetikzlibrary{arrows.meta}
\usetikzlibrary{calc}
\usetikzlibrary{decorations.pathmorphing}
\usetikzlibrary{external}
\usetikzlibrary{matrix}
\usetikzlibrary{patterns}
\usetikzlibrary{positioning}
\usetikzlibrary{quotes}
\usetikzlibrary{shapes}
\usetikzlibrary{shapes.misc}
\usetikzlibrary{spy}
\usepgfplotslibrary{groupplots}
\usepgfplotslibrary{fillbetween}
\pgfplotsset{colormap/viridis}
\pgfplotsset{compat=newest}
\tikzexternaldisable

\biboptions{sort&compress}
\AtBeginDocument{\hypersetup{
    linkcolor=[RGB]{0,117,162},
    citecolor=[RGB]{0,117,162},
    urlcolor=[RGB]{0,117,162}
}}

\definecolor{tol/contrast/blue}{RGB}{0,68,136}
\definecolor{tol/contrast/red}{RGB}{187,85,102}
\definecolor{tol/contrast/yellow}{RGB}{221,170,51}

\definecolor{tol/vibrant/blue}{RGB}{0,119,187}
\definecolor{tol/vibrant/cyan}{RGB}{51,187,238}
\definecolor{tol/vibrant/teal}{RGB}{0,153,136}
\definecolor{tol/vibrant/orange}{RGB}{238,119,51}
\definecolor{tol/vibrant/red}{RGB}{204,51,17}
\definecolor{tol/vibrant/magenta}{RGB}{238,51,119}
\definecolor{tol/vibrant/grey}{RGB}{187,187,187}

\definecolor{highlighter}{RGB}{222,244,64}

\DeclareCaptionSubType*{figure}

\newcommand{\latin}[1]{\textit{#1}}

\newcommand{\cmark}{\makebox[1em][c]{\textcolor{tol/vibrant/teal}{\ding{51}}}}
\newcommand{\xmark}{\makebox[1em][c]{\textcolor{tol/vibrant/red}{\ding{55}}}}
\newcommand{\midtilde}{\raisebox{-0.25\baselineskip}{\textasciitilde}}

\renewcommand{\vec}[1]{\boldsymbol{#1}}
\newcommand{\arr}[1]{\vv{#1}\vphantom{\ensuremath{#1}}}

\newcommand{\param}{\theta}
\newcommand{\q}[1][]{\if\relax\detokenize{#1}\relax\param\else\param^{#1}\fi}
\newcommand{\qs}{\vec{\param}}
\newcommand{\laplacian}{\Delta}

\newcommand{\grad}[1][]{\nabla\if\relax\detokenize{#1}\relax\else{}_{#1}\fi}
\newcommand{\timestep}{k}
\renewcommand{\div}[1][]{\grad[#1]\cdot}
\renewcommand{\d}{\relax\ifnum\lastnodetype>0\mskip\thinmuskip\fi\textnormal{d}}
\renewcommand{\L}{\mathcal{L}}
\DeclareMathOperator{\modulo}{mod}
\newcommand{\R}{\mathbb{R}}

\newcommand{\domain}{\ensuremath{\Omega}}
\newcommand{\interface}{\ensuremath{\Gamma}}

\newcommand{\sphere}{\ensuremath{\mathbb{S}^2}}
\newcommand{\Dirac}{\delta}
\newcommand{\Kronecker}{\delta}

\newcommand{\cardinality}{n}
\newcommand{\x}{\vec{x}}
\newcommand{\X}{\vec{X}}
\newcommand{\Xp}{\vec{\chi}}
\newcommand{\sites}{\ensuremath{\Theta}}
\newcommand{\metric}{g}
\newcommand{\weight}[1][]{\if\relax\detokenize{#1}\relax W\else W_{#1}}
\newcommand{\energy}{\mathcal{E}}
\newcommand{\kernel}{\hat{\delta}_1}

\newcommand{\density}{\rho}
\newcommand{\viscosity}{\mu}
\newcommand{\scinot}[2]{{#1}\times 10^{#2}}
\newcommand{\sci}[1]{\times10^{#1}}

\renewcommand{\u}{\vec{u}}
\newcommand{\floor}[1]{\left\lfloor{#1}\right\rfloor}
\newcommand{\U}{\vec{U}}
\newcommand{\f}{\vec{f}}
\newcommand{\F}{\vec{F}}
\newcommand{\e}{\vec{e}}
\newcommand{\n}{\vec{n}}

\newlength{\qendlen}
\newcommand{\qend}[1]{\setlength{\qendlen}{\widthof{#1}}#1\kern-1\qendlen}

\newcommand{\vanishingspace}{\relax\ifnum\lastnodetype>0\mskip\thinmuskip\fi}

\newcommand{\ns}{\vanishingspace\si{\nano\second}}
\newcommand{\um}{\vanishingspace\si{\micro\meter}}

\newcommand{\us}{\vanishingspace\si{\micro\second}}
\newcommand{\ms}{\vanishingspace\si{\milli\second}}

\newcommand{\umpersec}{\vanishingspace\si{\micro\meter\per\second}}
\newcommand{\mmpersec}{\vanishingspace\si{\milli\meter\per\second}}

\newcommand{\persec}{\vanishingspace\si{\per\second}}

\newcommand{\erg}{\vanishingspace\si{erg}}
\newcommand{\dynpercm}{\vanishingspace\si{dyn\per\centi\meter}}
\newcommand{\dynsecpercm}{\vanishingspace\si{dyn\cdot\second\per\centi\meter}}

\newcommand{\shear}{\dot{\gamma}}

\makeatletter
\newcommand{\macroname}[1]{\expandafter\@gobble\detokenize\expandafter{\string#1}}
\newcommand{\decorate}[2]{\expandafter\def\csname decor@\macroname{#1}\endcsname{#2}}
\newcommand{\decoration}[2]{\expandafter{#1}\csname decor@\macroname{#1}\endcsname{#2}}
\makeatother
\newcommand{\subscript}[1]{{}_{#1}}
\newcommand{\superscript}[1]{{}^{#1}}

\newcommand{\scsuperscript}[1]{\superscript{\textsc{#1}}}

\decorate{\sites}{\scsuperscript}
\decorate{\Xp}{\scsuperscript}
\decorate{\theta}{\scsuperscript}
\decorate{\varphi}{\scsuperscript}
\decorate{\cardinality}{\subscript}
\decorate{\vec{\psi}}{\scsuperscript}
\decorate{\vec{\bar{\psi}}}{\scsuperscript}
\decorate{\X}{\scsuperscript}

\newcommand{\data}[1]  {\decoration{#1}{d}}
\newcommand{\sample}[1]{\decoration{#1}{s}}
\newcommand{\poly}[1]  {\decoration{#1}{p}}
\newcommand{\reference}[1]{\hat{#1}}
\newcommand{\rbc}[1]{{#1}_\text{RBC}}
\newcommand{\plt}[1]{{#1}_\text{plt}}
\newcommand{\ethm}[1]{{#1}_\text{endo}}

\newcommand{\Xpsj}{\sample\Xp_{j}}
\newcommand{\Xsj}{\sample\X_{j}}

\begin{document}

\begin{frontmatter}

\title{Immersed boundary simulations of cell-cell interactions in whole blood}

\author[1]{Andrew Kassen\corref{corr}}  \ead{kassen@math.utah.edu}
\author[1]{Aaron Barrett}               \ead{barrett@math.utah.edu}
\author[2]{Varun Shankar}               \ead{shankar@cs.utah.edu}
\author[1]{Aaron L. Fogelson\fnref{3}}  \ead{fogelson@math.utah.edu}

\cortext[corr]{Corresponding author}
\address[1]{Department of Mathematics, University of Utah, Salt Lake City, UT 84112, USA}
\address[2]{School of Computing, University of Utah, Salt Lake City, UT 84112, USA}
\address[3]{Department of Bioengineering, University of Utah, Salt Lake City, UT 84112, USA}

\begin{abstract}
We present a new method for the geometric reconstruction of elastic surfaces simulated by the immersed boundary
method with the goal of simulating the motion and interactions of cells in whole blood. Our method uses
parameter-free radial basis functions for high-order meshless parametric reconstruction of point clouds and the
elastic force computations required by the immersed boundary method. This numerical framework allows us to
consider the effect of endothelial geometry and red blood cell motion on the motion of platelets. We find red
blood cells to be crucial for understanding the motion of platelets, to the point that the geometry of the vessel
wall has a negligible effect in the presence of RBCs. We describe certain interactions that force the platelets to
remain near the endothelium for extended periods, including a novel platelet motion that can be seen only in
3-dimensional simulations that we term ``unicycling\qend.'' We also observe red blood cell-mediated interactions
between platelets and the endothelium for which the platelet has reduced speed. We suggest that these behaviors
serve as mechanisms that allow platelets to better maintain vascular integrity.
\end{abstract}


\end{frontmatter}

\section{Introduction}

Blood is a complex mixture of cellular and fluid-phase components, composed most notably of red blood cells (RBCs)
and platelets suspended in plasma. RBCs are primarily the basis for transporting oxygen throughout the body.
Platelets, meanwhile, play a key role in the maintenance of the vasculature. Blood flows under pressure through
vessels, which vary in diameter between a few centimeters in the aorta to a few microns in the capillaries. A
single layer of endothelial cells, called the endothelium, lines healthy blood vessels.  Interactions between
platelets and RBCs, and between platelets and the endothelium, are important for a platelet's function, but models
of these interactions typically treat them in isolation.

At rest, RBCs are biconcave disk-shaped cells approximately $8\um$ in diameter and $2.5\um$ in thickness. The
volume fraction occupied by RBCs, or \emph{hematocrit}, ranges from approximately 36\% to 45\% in healthy humans.
In order to deliver oxygen throughout the body, RBCs must be extremely flexible, as some vessels are smaller than
the cell itself. Due to the wide variety of shapes exhibited by RBCs, their mechanical properties were studied
intensively during the 1970s and 80s. Canham~\cite{Canham:1970wx} theorized that the biconcave disk shape
minimizes bending energy. Skalak \latin{et al.}~\cite{Skalak:1973tp} devised a purpose-built constitutive law to
describe the tension of RBC membranes. Under the assumption of a viscoelastic response, Evans \& Hochmuth~%
\cite{Evans:1976tx} estimated the membrane viscosity. Mohandas \& Evans~\cite{Mohandas:1994tg} gave estimates of
the shear, bulk, and bending moduli, which have guided RBC models ever since~\cite{Pozrikidis:2003ft, Fai:2013do}.

Platelets in their inactive state are ellipsoidal disks, approximately 3--$4\um$ by $1\um$ in size. They are much
more rigid than RBCs due to their actin and microtubule-based cytoskeletons. Platelets are also much less
numerous, with 10--20 RBCs per platelet. Less is known about the mechanical properties of platelets. Models range
from perfectly rigid ellipsoids~\cite{Wang:2013gs} to systems of springs with~\cite{Erickson:2010ep,
Skorczewski:2013jn} or without~\cite{Wu:2014gt} a preferred curvature. One study estimates the shear modulus and
viscosity for platelets~\cite{Haga:1998wa}, but models tend to use a higher shear modulus than estimated and
neglect viscous effects altogether.

Because they are deformable, RBCs tend to move towards the center of a blood vessel, and in doing so may encounter
platelets, but the relative size and deformability of the RBC means a platelet is ejected from the RBC's path,
ultimately pushing the platelet into an RBC-free layer along the vessel wall. This process, called margination,
affects platelets and leukocytes (white blood cells) alike, and is the focus of many studies~\cite{Freund:2007kx,
Erickson:2010ep, Erickson:2011cf, Zhao:2011do, Kumar:2011dd, Zhao:2012ggba, Fedosov:2012dy, Kumar:2012ie,
Fedosov:2013ul, Muller:2014is, Fedosov:2014bs, Vahidkhah:2014hy, Vahidkhah:2015ch, Mehrabadi:2016fn}. From their
marginated positions, platelets survey the vessel wall for injury. Injury sites expose proteins, \latin{e.g.},
collagen and von Willebrand factor (vWF), to which platelets can bind and become activated. Platelet contact with
the injured wall is the essential first step.  This, in turn, leads the platelet to bind to the injury site and
release its own chemical signals to recruit further platelets, which eventually results in the formation of a
thrombus. All of this occurs in flowing blood, which sweeps these chemical signals downstream. While mechanisms
for platelet activation have been proposed for low and pathologically high shear rates, the case of
physiologically high shear rates is undecided~\cite{Fogelson:2015fb}.

Models of platelet motion over a thrombus indicate that there are stagnation zones immediately upstream and
downstream of the thrombus, where the fluid velocity is very slow, even when the thrombus protrudes only a few
microns from the vessel wall~\cite{Skorczewski:2013jn,Wang:2013gs}. Platelets that enter these regions may spend
an extended period of time near the thrombus. The portion of an endothelial cell containing its nucleus also
protrudes into the vessel approximately $1\um$. Moreover, these endothelial bumps are roughly periodic. If the
endothelium creates a stagnation zone, the trailing zone from one protrusion might lead into the leading zone of
the subsequent protrusion. This may allow for the sequestration of platelets or chemical signals. However, typical
models of platelet-wall interaction model the endothelium as a flat surface~\cite{Wu:2014gt,Vahidkhah:2015ch}.

The goal of this article is therefore to conduct 3D simulations of whole blood, incorporating red blood cell,
platelet, and endothelium interactions. Our model treats platelets and red blood cells as discrete elastic objects
immersed in and interacting with blood (which is modeled as an incompressible Newtonian fluid). We use this model
to compare the flow of whole blood across bumpy and flat walls and characterize the behavior and interactions of
platelets with the wall and RBCs. To simulate this model, we develop a cohesive numerical framework comprising a
fluid-structure interaction method, a meshless parametric modeling technique for reconstructing cell surfaces from
point clouds and computing elastic force densities on these surfaces, and a meshless quadrature scheme that
enables numerical integration of force densities on these surfaces as well.

We use the immersed boundary (IB) method for fluid-structure interaction. Originally developed by Peskin to study
the flow of blood around heart valves~\cite{Peskin:1972wa}, it has since been used to simulate, among numerous
other applications, vibrations in the inner ear~\cite{BeyerJr:1990tb}, the opening of a porous parachute~%
\cite{Kim:2006ku}, and sperm motility~\cite{Dillon:2011cu}, and has generated numerous related methods. The IB
method remains popular for modeling fluid-structure interaction because of its simplicity and ease of use, and
involves maintaining an Eulerian description of the fluid and a purely Lagrangian description of all immersed
elastic structures.

For parametric modeling and force density computations on these immersed elastic structures, we utilize meshless
interpolation based on radial basis functions (RBFs), which have been used for generating differentiation matrices
for the solution of PDEs~\cite{Fasshauer:2007ui}, surface reconstruction~\cite{Hardy:1971tb, Carr:2001tb,
Shankar:2013ki, SFKSISC2018}, and in the context of regularized Stokeslets to represent interfaces and approximate
their geometries~\cite{Olson:2015ja}. More relevantly, RBFs have been used in the context of the IB method to
reconstruct platelet surfaces from point clouds and to compute Lagrangian force densities in 2D simulations~%
\cite{Shankar:2015km}. However, this RBF-IB method has yet to be applied to surface reconstruction and force
density calculation in 3D simulations. Further, the 2D version of the RBF-IB method presented in~%
\cite{Shankar:2015km} required tuning in the RBF representation to achieve stability. In this work, we present the
first extension of the RBF-IB method to the simulation of whole blood in 3D geometries. Due to the use of the
meshless high-order accurate RBF-based representation, we are able to represent the constituent cells within whole
blood as point clouds with relatively small cardinality.  Further, we eliminate the aforementioned tuning
parameter using a recently-developed parameter-free RBF representation~\cite{SFKSISC2018}.

The remainder of this paper is organized as follows. We begin with an overview of the IB method in \cref{sec:ib}.
We then describe our method for solving the incompressible Navier-Stokes equations in \cref{sec:ins}. We describe
the elastic models used for each type of cell in \cref{sec:energy} and \cref{sec:rbfs} details our methods for
discretizing the cells using RBFs. Our results are presented in \cref{sec:results}.  Finally, we discuss the
implications of our findings in \cref{sec:conclusion}.

\section{The immersed boundary method}\label{sec:ib}

\subsection{Overview}\label{sec:ib_old}

Consider a rectangular domain, $\domain\subset\R^3$, which contains one or more deformable structures and is
otherwise filled with an incompressible, Newtonian fluid with constant density, $\density$, and viscosity,
$\viscosity$. The IB method treats these structures as an extension of the fluid. The motion of any particle in
$\domain$ is therefore governed by the incompressible Navier-Stokes equations,
\begin{gather}
    \density\left(\frac{\partial\u}{\partial t} + \div(\u\otimes\u)\right) = \viscosity\laplacian \u - \grad p + \f, \label{eq:ins-momentum} \\
    \div \u = 0, \label{eq:ins-incomp}
\end{gather}
where, for $\x = (x, y, z) \in \domain$, $\u = \u(\x, t) = (u, v, w)$ is the fluid velocity, $p = p(\x, t)$ is the
pressure, and $\f = \f(\x, t)$ is an external force density. Here and throughout this paper, we use bold italic
symbols to indicate vectors in $\R^3$. Treating the entire domain as a fluid allows us to discretize $\domain$
independently of the immersed structures with a fixed Eulerian grid of spacing $h$. The discretization of~%
\eqref{eq:ins-momentum} and~\eqref{eq:ins-incomp} is discussed in detail in \cref{sec:ins}.

Let $\X = \X(\theta, \varphi, t)$, for surface coordinates $(\theta, \varphi)$ in $\sites\subset\R^2$, be a
parametrization for the immersed boundary $\interface$. In the IB method, a Lagrangian representation is used to
track immersed boundaries.  Nevertheless, in continuum form, $\X$ satisfies a no-slip condition and moves with the
background fluid. The velocity of $\interface$ can therefore be represented as a convolution of the local fluid
velocity against the Dirac delta function $\Dirac(\x)$, \latin{i.e.},
\begin{equation}\label{eq:ib-interp}
    \frac{\partial \X}{\partial t} = \U(\X, t) 
        = \int\limits_{\domain} \u(\x, t) \Dirac(\x-\X) \d\x,
\end{equation}
where $\U$ is the restriction of $\u$ to $\interface$.
The process of transferring the fluid velocity to the immersed structure is termed ``interpolation\qend.'' 
As a boundary deforms, it generates a force density $\F = \F(\X, t)$, which it imparts onto
the fluid as $\f$ in~\eqref{eq:ins-momentum}. Similar to~\eqref{eq:ib-interp}, $\F$ is transferred to the fluid at
$\x$ via
\begin{equation}\label{eq:ib-spread}
        \f(\x, t)
        = \int\limits_{\interface} \F(\X, t)\Dirac(\x-\X) \d\X.
\end{equation}
This transfer of forces is called ``spreading\qend.'' The structure $\interface$ is typically represented in
Lagrangian form, usually as a discrete set of points. Consequently,~\eqref{eq:ib-interp} must be discretized to
update those points based on their velocities. Similarly,~\eqref{eq:ib-spread} must be discretized to both compute
$\F(\X, t)$ and to approximate the integral of $\F(\X, t)$ over $\interface$. Finally, the Dirac delta
$\Dirac(\x)$ is typically replaced by a smoothed, compactly-supported, $h$-dependent analog, typically referred to
as the ``discrete Delta function'' $\Dirac_h(\x)$.

In the classical version of the IB method~\cite{Peskin:2002go}, velocities are interpolated to the same Lagrangian
points that forces are spread from.  However, there exist several IB methods which instead use different sets of
points for each operation instead. For instance, Griffith \& Luo~\cite{Griffith:2017id} use a finite element
representation for the structure, and consequently spread forces from (Lagrangian) quadrature points, interpolate
velocities to quadrature points, and then project these using the finite element basis to the individual element
nodes (which are also Lagrangian points). The RBF-IB method~\cite{Shankar:2015km} used in this work does something
similar. In \cref{sec:rbfib}, we discuss the discretization of~\eqref{eq:ib-interp} and~\eqref{eq:ib-spread} in
the context of the RBF-IB method.

\subsection{The RBF-IB method}\label{sec:rbfib}

The RBF-IB method obtains its name from the fact that it uses meshless interpolation with RBFs to parametrize the
structure $\interface$ as $\X(\theta, \varphi, t)$. However, it also distinguishes itself from the classical IB
method by using two sets of Lagrangian points to represent the structure. One set of points, which we term
``movement points'' and label $\data\X_k$, $k=1, \ldots, \data\cardinality$, are used to move the structure in the
discrete analog to~\eqref{eq:ib-interp}. Thus, the fully discrete version of~\eqref{eq:ib-interp} in the RBF-IB
method reads
\begin{equation}\label{eq:ib-interp-disc}
    \frac{\partial\data\X_k}{\partial t} = \U(\data\X_k, t) \approx h^3 \sum_i \u(\x_i, t) \Dirac_h(\x_i-\data\X_k), \quad k=1, \ldots, \data\cardinality,
\end{equation}
where $h$ is the spacing of the background Eulerian grid and $i$ enumerates Eulerian grid points. 
In this work, we select the
movement points to be initially approximately $2h$ apart in $\R^3$ at the start of a simulation. The RBF-IB method
also uses a second set of points, which we term ``spreading points'' and label $\Xsj$,
$j=1, \ldots, \sample\cardinality$. These points are used to discretize~\eqref{eq:ib-spread} in order to spread
forces to the background Eulerian grid. The fully discrete version of~\eqref{eq:ib-spread} in the RBF-IB method
for an Eulerian grid point $\x_i$ can therefore be written as
\begin{align}\label{eq:ib-spread-disc}
        \f(\x_i, t) \approx\sum\limits_{j=1}^{\sample\cardinality} \weight[j]\F(\Xsj) \Dirac_h(\x_i-\Xsj),
\end{align}
where $\weight[j]$, $j=1, \ldots, \sample\cardinality$ is a set of quadrature weights for integrating quantities
at the spreading points $\sample\X$; the spreading points can therefore be thought of as a set of quadrature
points.

In general, the spreading points are selected to be approximately $h$ apart in $\R^3$ at the start of a
simulation; this spacing corresponds to spacing used for the single set of Lagrangian points (IB points) used in
traditional IB simulations. Thus, in the RBF-IB method, $\sample\cardinality > \data\cardinality$. In this sense,
the RBF-IB method is a generalization of the standard IB method. The use of movement points and spreading points
partially decouples the error in representing the structure $\interface$ from the error in integrating forces. The
spacing of the movement points $\data\X$ is chosen so that the reconstruction of $\interface$ from $\data\X$ and
the calculation of $\F(\X, t)$ by the RBF geometric model is sufficiently accurate. On the other hand, the spacing
of the spreading points $\sample\X$ is chosen so that the numerical integration of $\F(\X, t)$ is sufficiently
accurate, to prevent leaks, and to ensure that spread forces overlap sufficiently on the Eulerian grid. In this
work, we identify a set of fixed parametric locations $\data\sites = \{(\theta_k, \varphi_k)\}$,
$k=1, \ldots, \data\cardinality$ with the movement points $\data\X = \X(\data\sites, t)$, and call them ``data sites\qend.''
Analogously, we also maintain a set of fixed parametric points $\sample\sites = \{(\theta_j, \varphi_j)\}$,
$j=1, \ldots, \sample\cardinality$, that we identify with the spreading points $\sample\X$; these points are
called ``sample sites\qend.'' In general, we need only explicitly track the movement points $\data\X$, reconstruct
$\interface$ from $\data\X$ and the $\data\sites$ using our RBF geometric model in order to compute $\F(\X, t)$,
and regenerate the spreading points $\sample\X$ at the end of every step using the $\sample\sites$ and the RBF
geometric model.

Despite the advantages mentioned above, the choice of two sets of Lagrangian points potentially raises certain
concerns, which we address here. The first major concern relates to the force-spreading and velocity interpolation
operations. In traditional IB methods, the corresponding spreading and interpolation operators are formally
adjoint, leading to conservation of energy/power~\cite{Peskin:2002go}. However, in the RBF-IB method, these
operators are no longer adjoint, leading to concerns about numerical stability. Nevertheless, in previous
work~\cite{Shankar:2015km}, the authors demonstrated that the RBF-IB method dissipated energy in 2D simulations as
expected, despite the absence of formal adjointness. In this work, we show similar results for an RBC relaxation
problem in a full 3D simulation in \cref{sec:energy-est}. The second major concern relates to the wider $2h$
spacing of the movement points $\data\X$. If the points are spaced too widely apart, it is possible to obtain
unphysical behavior, with some parts of the structure responding more to the neighboring fluid than others.
However, both in~\cite{Shankar:2015km} and in this work, we observed that starting these points with a spacing of
$2h$ was sufficient to reproduce standard results for RBC simulations. It may be possible to alleviate any spacing
issues by periodically rearranging the movement points $\data\X$, but this approach was not necessary in this
work, and so we leave a full exploration of such Lagrangian rearrangement strategies for future work.

On a higher level,~\eqref{eq:ib-spread-disc} illustrates the need for three pieces of information for each
immersed structure: the points $\Xsj$ used to evaluate forces on the structure; a force density $\F(\Xsj)$ at each
of those points; and surface integration or quadrature weights $\weight[j]$ for those points. We detail the
RBF-based meshless geometric model used to obtain $\Xsj$, $\F(\Xsj)$, and $\weight[j]$ in \cref{sec:rbfs}, give
analytic expressions for energy densities used to derive $\F$ in \cref{sec:energy}, and give analytic expressions
for the force densities $\F$ in~\ref{sec:forces}.

\section{Solution of the incompressible Navier-Stokes equations}\label{sec:ins}

\subsection{Spatial Discretization}\label{sec:ns_space}

\begin{figure}[tb]
    \centering
    \begin{subfigure}{0.33\textwidth}
        \centering
    \begin{tikzpicture}[scale=1.5]
        \draw[help lines] (-0.3, -0.3) grid (2.3, 2.3);

        \foreach \x in {0,...,2}
            \foreach \y in {0,...,1}
            {
                \node[circle, inner sep=0pt, color=black!50, fill=white] at (\x, \y+0.5) {$u$};
            }
        \foreach \y in {0,...,2}
        {
            \foreach \x in {0,...,2}
                \node[circle, inner sep=0pt, fill=white] (av\x\y) at (\x, \y) {$\bar{v}$};
            \foreach \x in {0,...,1}
                \node[circle, inner sep=0pt, fill=white] (v\x\y) at (\x+0.5, \y) {$v$};
        }
        \path[-stealth] (v00.south west) edge[bend left] node[midway, below, yshift=+1pt] {} (av00.south east);
        \path[-stealth] (v10.south west) edge[bend left] node[midway, below, yshift=+1pt] {} (av10.south east);
        \path[-stealth] (v01.south west) edge[bend left] node[midway, below, yshift=+1pt] {} (av01.south east);
        \path[-stealth] (v11.south west) edge[bend left] node[midway, below, yshift=+1pt] {} (av11.south east);
        \path[-stealth] (v02.south west) edge[bend left] node[midway, below, yshift=+1pt] {} (av02.south east);
        \path[-stealth] (v12.south west) edge[bend left] node[midway, below, yshift=+1pt] {} (av12.south east);

        \path[-stealth] (v00.north east) edge[bend left] node[midway, below, yshift=+1pt] {} (av10.north west);
        \path[-stealth] (v10.north east) edge[bend left] node[midway, below, yshift=+1pt] {} (av20.north west);
        \path[-stealth] (v01.north east) edge[bend left] node[midway, below, yshift=+1pt] {} (av11.north west);
        \path[-stealth] (v11.north east) edge[bend left] node[midway, below, yshift=+1pt] {} (av21.north west);
        \path[-stealth] (v02.north east) edge[bend left] node[midway, below, yshift=+1pt] {} (av12.north west);
        \path[-stealth] (v12.north east) edge[bend left] node[midway, below, yshift=+1pt] {} (av22.north west);

        \node[black] at (-0.25, 2.25) {(a)};
    \end{tikzpicture}
    \label{fig:adv-x-ave}
    \end{subfigure}%
    \begin{subfigure}{0.33\textwidth}
        \centering
    \begin{tikzpicture}[scale=1.5]
        \draw[help lines] (-0.3, -0.3) grid (2.3, 2.3);

        \foreach \x in {0,...,2}
        {
            \foreach \y in {0,...,1}
            {
                \node[circle, inner sep=0pt, fill=white] (u\x\y) at (\x, \y+0.5) {$u$};
            }
            \foreach \y in {0,...,2}
            {
                \node[circle, inner sep=0pt, fill=white] (au\x\y) at (\x, \y) {$\bar{u}$};
            }
        }
        \foreach \x in {0,...,1}
            \foreach \y in {0,...,2}
            {
                \node[circle, inner sep=0pt, color=black!50, fill=white] at (\x+0.5, \y) {$v$};
            }

        \path[-stealth] (u00.south east) edge[bend left] node[midway, below, yshift=+1pt] {} (au00.north east);
        \path[-stealth] (u10.south east) edge[bend left] node[midway, below, yshift=+1pt] {} (au10.north east);
        \path[-stealth] (u20.south east) edge[bend left] node[midway, below, yshift=+1pt] {} (au20.north east);
        \path[-stealth] (u01.south east) edge[bend left] node[midway, below, yshift=+1pt] {} (au01.north east);
        \path[-stealth] (u11.south east) edge[bend left] node[midway, below, yshift=+1pt] {} (au11.north east);
        \path[-stealth] (u21.south east) edge[bend left] node[midway, below, yshift=+1pt] {} (au21.north east);

        \path[-stealth] (u00.north west) edge[bend left] node[midway, below, yshift=+1pt] {} (au01.south west);
        \path[-stealth] (u10.north west) edge[bend left] node[midway, below, yshift=+1pt] {} (au11.south west);
        \path[-stealth] (u20.north west) edge[bend left] node[midway, below, yshift=+1pt] {} (au21.south west);
        \path[-stealth] (u01.north west) edge[bend left] node[midway, below, yshift=+1pt] {} (au02.south west);
        \path[-stealth] (u11.north west) edge[bend left] node[midway, below, yshift=+1pt] {} (au12.south west);
        \path[-stealth] (u21.north west) edge[bend left] node[midway, below, yshift=+1pt] {} (au22.south west);
        \node[black] at (-0.25, 2.25) {(b)};
    \end{tikzpicture}
    \label{fig:adv-y-ave}
    \end{subfigure}%
    \begin{subfigure}{0.33\textwidth}
        \centering
    \begin{tikzpicture}[scale=1.5]
        \draw[help lines] (-0.3, -0.3) grid (2.3, 2.3);

        \foreach \x in {0,...,2}
            \foreach \y in {0,...,2}
                \node[circle, inner sep=0pt, fill=white] (uv\x\y) at (\x, \y) {$\bar{u}\bar{v}$};
        \foreach \x in {0,...,2}
            \foreach \y in {0,...,1}
                \node[circle, inner sep=0pt] (duv\x\y) at (\x, \y+0.5) {$\star$};

        \path[stealth-] (duv00.south east) edge[bend left] (uv00.north east);
        \path[stealth-] (duv10.south east) edge[bend left] (uv10.north east);
        \path[stealth-] (duv20.south east) edge[bend left] (uv20.north east);
        \path[stealth-] (duv01.south east) edge[bend left] (uv01.north east);
        \path[stealth-] (duv11.south east) edge[bend left] (uv11.north east);
        \path[stealth-] (duv21.south east) edge[bend left] (uv21.north east);

        \path[stealth-] (duv00.north west) edge[bend left] (uv01.south west);
        \path[stealth-] (duv10.north west) edge[bend left] (uv11.south west);
        \path[stealth-] (duv20.north west) edge[bend left] (uv21.south west);
        \path[stealth-] (duv01.north west) edge[bend left] (uv02.south west);
        \path[stealth-] (duv11.north west) edge[bend left] (uv12.south west);
        \path[stealth-] (duv21.north west) edge[bend left] (uv22.south west);
        \node[black] at (-0.25, 2.25) {(c)};
    \end{tikzpicture}
    \label{fig:adv-y-diff}
    \end{subfigure}%
    \caption{%
A cross-section illustrating the steps in computing the $D_2[(A_2 u) (A_1 v)]$ term of the first component of the
advection. The horizontal and vertical velocity component discretization locations are marked by $u$ and $v$,
respectively. Arrows emanate from a point contributing to a stencil and point to the center of the stencil.
(a) $A_1$ averages $v$ in the $x$ direction, yielding an approximation $\bar{v}$ at grid vertices (in 3D, centers
of cell edges for which $x$ and $y$ are constant).
(b) $A_2$ averages $u$ in the $y$ direction, yielding an approximation $\bar{u}$ at the same points as (a). The
quantities $A_1 v$ and $A_2 u$ are collocated and can be directly multiplied to obtain an approximation of $uv$ at
locations marked $\bar{u}\bar{v}$.
(c) $D_2$ approximately differentiates $uv$ in the $y$ direction, yielding the desired quantity at each point
marked $\star$.  The approximation of $uv$ is also used to compute $D_1[(A_1 v)(A_2 u)]$ in the second component
of the advection, wherein application of $D_1$ instead yields approximations collocated with locations marked $v$
in (a).
    }
    \label{fig:discretization}
\end{figure}
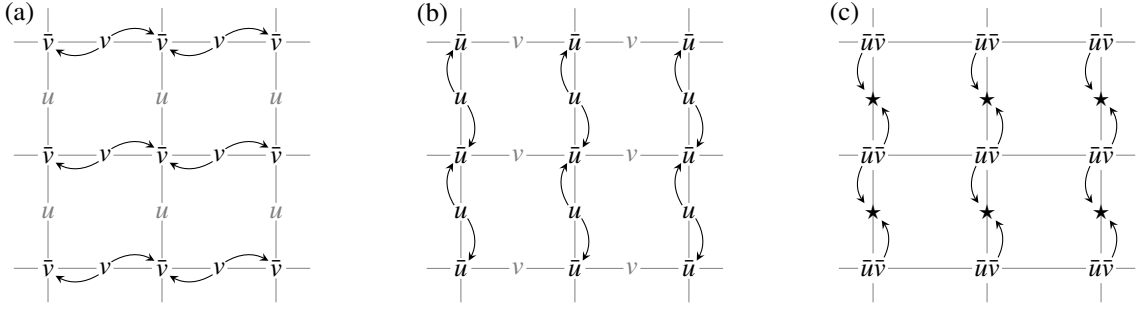

To discretize the Navier-Stokes equations~\eqref{eq:ins-momentum} and~\eqref{eq:ins-incomp}, we use a
marker-and-cell (MAC) grid~\cite{Welch:1965jv}: for grid cell center $\x_i$, scalar-valued function $s(\x)$ is
discretized at $\x_i$, and component $\e_a\cdot\vec{v}$ of vector-valued function $\vec{v}(\x)$ at
$\x_i-\sfrac12h\e_a$, where $\e_a$ is a canonical basis vector. Define the centered difference operator
\begin{equation}
    D_a\phi(\x) = \frac{\phi(\x+\sfrac12 h\e_a) - \phi(\x-\sfrac12 h\e_a)}{h},
\end{equation}
for which, \latin{e.g.}, $D_1$ approximates differentiation in the $x$ direction. The discrete divergence,
gradient, and Laplacian operators use centered differences, resulting in a 2-point stencil for each discrete first
derivative and the standard 7-point discrete Laplacian. We also define the centered average operator
\begin{equation}
    A_a\phi(\x) = \frac{\phi(\x+\sfrac12 h\e_a) + \phi(\x-\sfrac12 h\e_a)}2.
\end{equation}
By averaging $u$ in the $y$ direction and $v$ in the $x$ direction, we obtain collocated approximations to $u$ and
$v$ at the center of a cell edge. Averaging, e.g., $u$ in the $x$ direction yields an approximation to $u$ at the
cell center. We can therefore discretize the components of the advection term $\div(\u\otimes\u)$ by
\begin{equation}\label{eq:advection}
    \div[h](\u\otimes\u) :=
    \begin{bmatrix}
        D_1[(A_1 u) (A_1 u)] + D_2[(A_1 v) (A_2 u)] + D_3[(A_1 w) (A_3 u)] \\
        D_1[(A_2 u) (A_1 v)] + D_2[(A_2 v) (A_2 v)] + D_3[(A_2 w) (A_3 v)] \\
        D_1[(A_3 u) (A_1 w)] + D_2[(A_3 v) (A_2 w)] + D_3[(A_3 w) (A_3 w)]
    \end{bmatrix}.
\end{equation}
The symbol $\div[h]$ represents the discrete divergence operator.  \Cref{fig:discretization} illustrates the steps
in computing $D_2[(A_1 v)(A_2 u)]$, which appears in the first component in~\eqref{eq:advection}. Morinishi
\latin{et al.}~\cite{Morinishi:1998us} show that this scheme, $Div. - S2$ in their parlance, is conservative under
the assumption that $\u$ is discretely divergence-free.

\subsection{Temporal Discretization}\label{sec:ns_time}

To advance the solution, we use either the backward-forward Euler-based scheme~\cite{Ascher:1997tm} or the 2-stage
scheme described by Peskin~\cite{Peskin:2002go}, modified to advance structures using the newest velocities. Both schemes utilize a pressure projection (or fractional step) method. The
modification makes these schemes formally first-order in time, but allow us to separate the Eulerian update from
the Lagrangian update by requiring only information at the beginning of the timestep to evaluate forces and moving
the structure at the end of the timestep. For the backward-forward Euler scheme, discretizing~%
\eqref{eq:ins-momentum} to advance time to $t+\timestep$, yields linear solves of Helmholtz type,
\begin{equation}\label{eq:disc-momentum}
    (I - \timestep\rho^{-1}\mu \laplacian_h) \u^\ast = \u^n - \timestep\left[\div[h]\left(\u^n\otimes\u^n\right) + \rho^{-1}\left(\f^{n+1} - \grad[h]p^n\right)\right] \quad\text{in}~\domain, \\
\end{equation}
with boundary conditions
\begin{equation}\label{eq:disc-bdy}
    \u^\ast = \u^{n+1}_b + \timestep \grad[h] q^{n} \quad\text{on}~\partial\domain,
\end{equation}
where superscripts denote the time step, $\u_b$ is velocity boundary data, $\laplacian_h$ and $\grad[h]$ are the
discrete Laplacian and gradient, respectively, and $q$ is described below. The force density $\f^{n+1}$ is
advanced using only data at timestep $n$. The intermediate velocity field $\u^\ast$ may not be divergence-free. To
obtain a velocity field that satisfies~\eqref{eq:ins-incomp}, we use projection method II (PmII) of Brown, Cortez,
and Minion~\cite{Brown:2001bq}. PmII updates the pressure using
\begin{equation}
    p^{n+1} = p^n + (\rho I - \timestep\mu\laplacian_h)q^{n+1},
\end{equation}
and generates the divergence-free velocity field
\begin{equation}\label{eq:vel-update}
    \u^{n+1} = \u^\ast - \timestep\grad[h]q^{n+1}
\end{equation}
using the pseudo-pressure $q^{n+1}$, which satisfies
\begin{equation}
\begin{alignedat}{2}
    &k \laplacian_h q^{n+1} = \div[h]\u^\ast &&\quad \text{in}~\domain, \\
    &\n\cdot\grad[h] q^{n+1} = 0                    &&\quad \text{on}~\partial\domain.
\end{alignedat}
\end{equation}
The velocity update~\eqref{eq:vel-update} provides the boundary conditions~\eqref{eq:disc-bdy} using a lagged
value of the pseudo-pressure. The 2-stage RK method consists of a backward-forward Euler step followed by a
Crank-Nicolson-midpoint step, which involves only minor modifications to~\eqref{eq:disc-momentum}. In total, we
perform 3 Helmholtz solves and 1 Poisson solve per RK stage.

We employ preconditioned conjugate gradients (PCG) to perform the solves. We use Chebyshev iteration as a
preconditioner for the Helmholtz solves and as an error smoothing procedure and direct solver for multigrid (MG)
to precondition the Poisson solve. Chebyshev iteration is a generalization of weighted Jacobi iteration which
requires only the ability to perform sparse matrix polynomial-vector multiplication. Chebyshev iteration (MG) PCG
is therefore parallelized by using a parallel sparse matrix-vector multiplication routine with Horner's method to
evaluate the polynomials.

In the case of a triply periodic domain, the linear solves involve symmetric matrices. For Dirichlet or Neumann
boundaries, we extrapolate using values at the neighboring grid points and boundary data to fill ghost points. In
these situations, the standard discrete second derivative may actually approximate some non-unit multiple of its
continuous counterpart near the boundary. To account for this while maintaining symmetry of the Helmholtz
matrices, we scale equations involving near-boundary values. For details, see~\ref{sec:boundary-correction}. The
trade-off is 3 extra diagonal matrix-vector multiplications per RK stage for the ability to use PCG for the linear
solves.

\section{Cell energy and power models}\label{sec:energy}

In this section, we describe the various forms of energy density (energy per area) and power density (power per
area) used in our simulations, and give analytic expressions for each. The corresponding force densities are given
in~\ref{sec:forces}. We use five kinds of densities: spring energy, damped spring power, tension energy,
dissipative power, and Canham-Helfrich bending energy. For constitutive law $W$, we define the functional
\begin{equation}
    \energy[\X, \U] = \int\limits_\interface W(\X, \U, \ldots) \d\X,
\end{equation}
where the ellipsis indicates that $W$ may depend on spatial derivatives of $\X$ or $\U$. The force density
associated with $W$ is found by computing the first variation of $\energy$,
\begin{equation}
    \F = -\delta\energy,
\end{equation}
with respect to $\X$ for energy densities and to $\U$ for power densities. Because our ultimate goal is a
three-dimensional simulation, we limit our descriptions to the three-dimensional case. Considerations for the
two-dimensional case are treated elsewhere~\cite{Peskin:2002go,Erickson:2010uzba}.

We begin with Hookean energy and damped spring power density. These have the simplest constitutive laws we use,
depend only on surface locations and surface velocities, and take the form
\begin{equation}
        W_\text{Hk}(\X) = \frac{k}2 {\|\X - \X'\|}^2\quad\text{and}\quad
        W_\text{damped}(\U) = \frac{\eta}2{\|\U - \U'\|}^2,
\end{equation}
where $\X'=\X'(\theta, \varphi, t)$ is the tether location for $\X(\theta, \varphi, t)$,
$\U'=\U'(\theta, \varphi, t)$ is the prescribed velocity of the tether point, $k$ is the spring constant, and
$\eta$ is the damping constant.  Due to the lack of information about the mechanical properties of endothelial
cells, we model the endothelium as a rigid, stationary object with $\ethm{k}=2.5\dynpercm$ and
$\ethm\eta=2.5\sci{-7}\dynsecpercm$, chosen to be as large as possible for the chosen spatial and temporal step
size with prescribed velocity $\U' = \vec{0}$. We compare different choices for $\X'$ in \cref{sec:whole-blood}.

Next, we consider the tension energy densities for RBCs and platelets. These penalize stretching and areal
dilation of the cell membranes. Let $\lambda_1$ and $\lambda_2$ be the principal extensions, \latin{i.e.}, the
maximal and minimal ratios of stretching relative to a reference configuration. We define the invariants
$I_1=\lambda_1^2+\lambda_2^2-2$ and $I_2 = \lambda_1^2\lambda_2^2-1$, which measure relative changes in length and
area, respectively, such that $I_1 = I_2 = 0$ correspond to a rigid body motion. We express the tension density in
terms of these invariants. Skalak's Law was designed specifically for RBCs~\cite{Skalak:1973tp}:
\begin{equation}\label{eq:skalak-law}
    W_\text{Sk}(I_1, I_2) = \frac{E}4\left(I_1^2 + 2I_1 - 2I_2\right) + \frac{G}4 I_2^2.
\end{equation}
$E$ is the shear modulus, and $G$ is the bulk modulus. For RBCs, we follow Fai \latin {et al.}~\cite{Fai:2013do} and set $\rbc{E} = 2.5\sci{-3}\dynpercm$ and $\rbc{G} = 2.5\sci{-1}\dynpercm$. We use
the shape given by Evans \& Fung~\cite{Evans:1972uf} for the reference RBC with radius $R_0 = 3.91\um$,
\begin{equation}
    \rbc{\vec{\hat{X}}}(\theta, \varphi) = R_0\begin{bmatrix}
            \cos\theta\cos\varphi \\
            \sin\theta\cos\varphi \\
            z(\cos\varphi)\sin\varphi
    \end{bmatrix},
\end{equation}
where $(\theta, \varphi)\in(-\pi, \pi]\times[-\pi/2, \pi/2]$ and $z(r) = 0.105 + r^2 - 0.56r^4$. Platelets, on the
other hand, do not have a purpose-built constitutive law, but are known to be stiffer than RBCs. We use the neo-Hookean model
\begin{equation}\label{eq:neohookean}
    W_\text{nH}(I_1, I_2) = \frac{E}2\left(\frac{I_1+2}{\sqrt{I_2+1}}-2\right) + \frac{G}2 {\left(\sqrt{I_2+1}-1\right)}^2
\end{equation}
with $\plt{E} = 1\sci{-1}\dynpercm$ and $\plt{G} = 1\dynpercm$, and an ellipsoidal reference configuration~%
\cite{Frojmovic:1982wk}
\begin{equation}
    \plt{\hat{\vec{X}}}(\theta, \varphi) = \begin{bmatrix}
            1.55\um\cos\theta\cos\varphi \\
            1.55\um\sin\theta\cos\varphi \\
            0.5\um\sin\varphi
    \end{bmatrix}.
\end{equation}

Platelets and RBCs also respond to changes in membrane curvature. Let $H$ be the membrane's mean curvature. The
Canham-Helfrich bending energy density takes the form~\cite{Canham:1970wx}
\begin{equation}\label{eq:bending-energy}
    W_\text{CH}(H) = 2\kappa {(H-H')}^2,
\end{equation}
where $\kappa$ is the bending modulus in units of energy, and $H'$ is the spontaneous or \emph{preferred}
curvature. An RBC generates a relatively weak response to changes in its curvature. Its bending modulus is
estimated to be in the range $0.3$--$4\sci{-12}\erg$~\cite{Mohandas:1994tg}. We use a bending modulus of
$\rbc\kappa = 2\sci{-12}\erg$ and a preferred curvature $H' = 0$ for RBCs. RBCs, therefore, tend to locally
flatten their membranes. For platelets, we use a larger bending modulus of $\plt\kappa=2\sci{-11}\erg$
and a preference for its reference curvature. Together with the neo-Hookean tension above, this maintains a fairly
rigid platelet.

Finally, we consider dissipative power, which causes the membrane to exhibit a viscoelastic response to strain. It
takes the form~\cite{Rangamani:2012hi}
\begin{equation}\label{eq:dissip-energy}
    W_\text{dissip}(\dot{\lambda}_1, \dot{\lambda}_2) = \frac{\nu}{2}\left(\frac{\dot{\lambda}_1^2}{\lambda_1^2} + \frac{\dot{\lambda}_2^2}{\lambda_2^2}\right),
\end{equation}
where $\nu$ is the membrane viscosity, and $\dot{\lambda}_i$ is the rate of change of $\lambda_i$. We imbue only
the RBC with viscoelasticity. We find this effective in eliminating some numerical instabilities. While Evans \&
Hochmuth suggest a viscosity of approximately $1\sci{-3}\dynsecpercm$~\cite{Evans:1976tx}, we find this to be
prohibitively expensive in practice, due to time step restrictions, and instead use
$\rbc\nu = 2.5\sci{-7}\dynsecpercm$.

\section{Geometry of reconstructed surfaces}\label{sec:rbfs}

From the previous section, we have analytic expressions for the energy and power densities, and consequently the
Lagrangian force density, $\F$ (see~\ref{sec:forces}).  The RBC and platelet force models require first and second
derivatives along their surfaces.  Additionally, the IB force spreading operation requires quadrature weights for
the cell surfaces. This section finishes where we left off by discussing the construction of the necessary
discrete linear operators and quadrature weights $\weight[j]$ through the use of RBF-based methods.

\subsection{Surface reconstruction with radial basis functions}\label{sec:rbf-interpolation}

We now describe our RBF method for reconstructing cell surfaces. RBF interpolation is a meshfree approach to
scattered data interpolation where structural information is encoded purely as point-wise distances. This
contrasts with, \latin{e.g.}, polynomials, where points must be chosen at grid vertices, or spherical harmonics,
for which special node sets are typically used. RBFs have been shown to be a viable approach for representing
cells on a par with Fourier methods~\cite{Shankar:2013ki}. They are therefore appealing for representing red blood
cells and platelets.

Our goal is to produce a parametric reconstruction of cells. We choose to parametrize both red blood cells and
platelets on the 2-sphere, $\sphere$. We use a spherical coordinate mapping to relate Cartesian points on the
2-sphere to their corresponding parametric points. This mapping is given by
\begin{equation}
    \Xp(\theta, \varphi) =
    \begin{bmatrix}
        \cos\theta\cos\varphi \\
        \sin\theta\cos\varphi \\
        \sin\varphi
    \end{bmatrix},
\end{equation}
$(\theta, \varphi)\in(-\pi, \pi]\times[-\pi/2, \pi/2]$. As mentioned in Section~%
\ref{sec:rbfib}, let $\data\sites = \{(\theta_k, \varphi_k)\}$,
$k=1, \ldots, \data\cardinality$, be a set of distinct \emph{data sites}, defined by the
Bauer spiral~\cite{Bauer:2000km},
\begin{equation}\label{eq:bauer-spiral}
    \begin{aligned}
        &\varphi_k = \sin^{-1}(-1 + (2k - 1) / \data\cardinality),\\
        &\theta_k = \modulo\left(\sqrt{\data\cardinality\pi}\varphi_k + \pi, 2\pi\right) - \pi,
    \end{aligned}
\end{equation}
where $\modulo(a, b) = a - b\floor{a/b}$ is the modulo function. These data sites correspond to the movement points $\data\X$ to which the fluid grid velocity is interpolated. We similarly generate a larger set of
\emph{sample sites} from which to spread forces by replacing $\data\cardinality$ with $\sample\cardinality$ in the
above equation. Let $\data\Xp_i=\Xp(\theta_i, \varphi_i)$ for each $(\theta_i, \varphi_i)\in\data\sites$. In this
setting, our goal is to construct a parametric mapping from $\data\sites$ to the Cartesian locations of points on
cell surfaces. For RBC and platelet parametrizations, we identify the point $\X(\theta, \varphi, t)$ on
$\interface$ with the point $\Xp(\theta, \varphi)$ on $\sphere$. Each component of $\X$ and consequently each
component of a movement point $\data\X_k$ is then a function defined on $\sphere$. The problem of surface
reconstruction therefore involves approximating each of these functions from their values at $\data\sites$ using an RBF
interpolant. For the discussion that follows, we use $\psi(\Xp):\sphere\to\R$ to denote a function that we wish to
approximate.

Let $\phi:\sphere\times\sphere\to\R$ be a \emph{radial kernel} with the property that
$\phi(\Xp_i, \Xp_j)\equiv\phi(\|\Xp_i-\Xp_j\|)$. These kernels are sometimes called \emph{spherical basis
functions}. In addition to the radial kernels, let $p_k(\Xp)$, $k=1, \ldots, \poly\cardinality$ denote the first
$\poly\cardinality$ spherical harmonics, which form a natural basis for polynomial approximation on $\sphere$.
Then, the RBF interpolant to $\psi(\Xp)$ takes the form
\begin{equation}\label{eq:rbf-interp}
    s(\Xp)
    = \sum_{k=1}^{\data\cardinality} c_k \phi(\|\Xp-\data\Xp_k\|)
    + \sum_{i=1}^{\poly\cardinality} d_i p_i(\Xp),
\end{equation}
where $c_k$ and $d_i$ are unknown interpolation coefficients, and the metric can be written parametrically as
\begin{equation} \label{eq:param-metric}
    \|\Xp-\data\Xp_k\|
        = \|\Xp(\theta, \varphi) - \Xp(\theta_k, \varphi_k)\|
        = \sqrt{2(1 - \cos\varphi\cos\varphi_k\cos(\theta-\theta_k) - \sin\varphi\sin\varphi_k)}.
\end{equation}
To find $c_k$ and $d_i$, we enforce that the interpolant~\eqref{eq:rbf-interp} exactly interpolate the function
$\psi(\Xp)$ at each $\data\Xp_k$,
\begin{equation}
    s(\data\Xp_k) = \psi(\data\Xp_k), \qquad k=1, \ldots, \data\cardinality,\label{eq:interp_constraint}
\end{equation}
and that~\eqref{eq:rbf-interp} exactly reproduce the first $\poly\cardinality$ spherical harmonics everywhere on
$\sphere$~\cite{Fasshauer:2007ui},
\begin{equation}
    \sum_{k=1}^{\data\cardinality} c_k p_i(\data\Xp_k) = 0, \qquad i=1, \ldots, \poly\cardinality.\label{eq:constraints}
\end{equation}
We collect~\eqref{eq:interp_constraint} and~\eqref{eq:constraints} into a dense symmetric block system of the form
\begin{equation}\label{eq:rbf-interp-matrix}
    \begin{bmatrix}
            \Phi & P \\ P^T & 0
    \end{bmatrix}\begin{bmatrix}
            \arr{c} \\ \arr{d}
    \end{bmatrix} = \begin{bmatrix}
            \arr{\psi} \\ \arr{0}
    \end{bmatrix},
\end{equation}
where $\arr{c}$ and $\arr{d}$ are the unknown coefficients, $\arr{\psi}$ is the vector of evaluations of $\psi(\data\Xp)$, $\Phi$ represents the evaluations of $\phi$, $P$
represents evaluations of the polynomials $p_k$, the matrix block $0$ is the
$\poly\cardinality\times\poly\cardinality$ zero matrix, and $\arr{0}$ is a vector of $\poly\cardinality$ zeros.
Because $\data\sites$ is fixed, we construct this matrix and compute its factors only once even though $\psi(\data\Xp)$
changes. The matrix in~\eqref{eq:rbf-interp-matrix} is invertible for any conditionally positive definite kernel
$\phi$ of order $m$ as long as the data sites $\data\sites$ are unisolvent for the first
$\poly\cardinality \ge (m+1)^2$ spherical harmonics, \latin{i.e.}, $P$ is of full rank~ \cite{Fasshauer:2007ui}. A
common heuristic choice to ensure this is to set $\data\cardinality = 2 \poly\cardinality$~\cite{SWJCP2018}.

It now remains to discuss the choice of $\phi$. In previous work~\cite{Shankar:2015km}, we chose $\phi$ to be an
infinitely-smooth and positive-definite kernel. While these kernels offer spectral convergence rates, they require
tuning a so-called shape parameter~\cite{Fasshauer:2007ui}. In particular, we instead choose $\phi$ to be a
polyharmonic spline, which is a conditionally-positive definite kernel with finite smoothness. In this work,
we set $\phi(r) = r^7$, augmented with fifth-order spherical harmonics to represent RBCs, or just the constant
spherical harmonic for platelets. As mentioned previously, we identify the point $\X(\theta, \varphi, t)$ on
$\interface$ with the point $\Xp(\theta, \varphi)$ on $\sphere$.  To reconstruct the surface $\X$ from the
movement points $\data\X$, we simply use the fact that each component of a particular movement point $\data\X_k$
is one of $\data\cardinality$ samples of a function of the form $\psi(\Xp)$. For instance, to reconstruct the
first component of $\X$, we set $\arr{\psi}$ to be the vector of $x$-coordinates of the movement
points $\data\X$, and find a set of coefficients corresponding to the $x$-component of $\X$.  This process is
repeated component-wise for the movement points, and is essentially a component-wise interpolation of the movement
points~\cite{Shankar:2015km}.

Clearly, by replacing the quantity $\arr{\psi}$ with samples of any function defined on the sphere, one can also
use~\eqref{eq:rbf-interp} to interpolate other quantities on the surface $\interface$, since all such quantities
are functions of the form $\psi(\Xp)$. In this work, we also use~\eqref{eq:rbf-interp} to compute force densities
required by the IB method. It is clear from~\eqref{eq:skalak-law}--\eqref{eq:dissip-energy} that computing the
force densities requires values of $I_1$, $I_2$, and $H$, among others.  These values are derived from the first
and second derivatives of $\X$, which can be obtained by analytically differentiating the RBF and spherical
harmonic bases.

In practice, we do not compute any coefficients in our simulations. For efficiency, it is possible to reformulate
the process of interpolation followed by differentiation as a single application of a discrete (nodal)
differential operator (or a differentiation matrix). We discuss this in the next section.

\subsection{Discrete linear surface operators}

Let $\L$ be a linear operator. In particular, we are interested in the first- and second-order partial
differential operators, $\partial/\partial\theta$, $\partial^2/\partial\theta\partial\varphi$, \latin{etc}. First,
note that since $s(\Xp) \approx \psi(\Xp)$, it is also true that $\L s(\Xp) \approx \L \psi(\Xp)$. Thus, in the
RBF-IB method, it is of interest to be able to apply $\L$ to $s(\Xp)$ and evaluate the resulting
function efficiently. To do so, we formulate all applications and evaluations of $\L s$ in terms of a
matrix-vector multiplication with the quantity $\arr{\psi}$. To see how this is done, first note that if we wish
to approximate samples of $\L s$ at a given set of sample sites $\sample\sites$, we can write
\begin{align}
\left.\L s (\Xp) \right|_{\sample\sites}
    = \sum_{k=1}^{\data\cardinality} c_k \L \phi(\|\Xp-\data\Xp_k\|)
    + \sum_{i=1}^{\poly\cardinality} d_i \L p_i(\Xp),
\end{align}
which can be expressed more compactly in terms of matrix-vector multiplications as
\begin{align}
\left.\L s (\Xp) \right|_{\sample\sites} =
    \begin{bmatrix}
    \L \Phi & \L P
    \end{bmatrix}\begin{bmatrix}
    \arr{c} \\
    \arr{d}
    \end{bmatrix}.
\end{align}
However, we can use~\eqref{eq:rbf-interp-matrix} to rewrite this as
\begin{align}
\left.\L s (\Xp) \right|_{\sample\sites} =
    \begin{bmatrix}
    \L \Phi & \L P
    \end{bmatrix}\begin{bmatrix}
    \Phi & P \\
    P^T & 0
    \end{bmatrix}^{-1}\begin{bmatrix}
    \arr{\psi} \\
    \arr{0}
    \end{bmatrix} =
    \begin{bmatrix}
    L & \ast
    \end{bmatrix}\begin{bmatrix}
    \arr{\psi}\\
    \arr{0}
    \end{bmatrix},
\end{align}
where $\L\Phi$ and $\L P$ represent evaluations of $\L\phi$ and $\L p_k$, respectively, and the
$\sample\cardinality \times \data\cardinality$ (dense) matrix $L$ is the discrete analog of $\L$. It is important
to note that $L$ is completely independent of the function $\psi$; in fact, it depends only on the functions
$\phi$ and $p_k$, the fixed data sites $\data\sites$, and the fixed sample sites $\sample\sites$. The block marked
by $\ast$ is multiplied by zeros, and can be discarded. We compute a separate $L$ for each operator $\L$ as a
preprocessing step, and simply apply these matrices to any function that needs to be evaluated or differentiated.
For instance, setting $\L$ to be the point evaluation operator and setting $\arr{\psi}$ to each component of the
movement points $\data\X$ allows us to generate each component of the spreading points $\sample\X$.  Similarly, we
can obtain samples of parametric derivatives of any function $\psi(\Xp)$ at $\sample\sites$ by replacing $\L$ with
that derivative operator.

It is also straightforward to generate versions of $L$ that produce derivatives at $\data\sites$ simply by
replacing $\left.\L s\right|_{\sample\sites}$ with $\left.\L s\right|_{\data\sites}$ in the above discussion. With
these operators in hand, the quantities $I_1$, $I_2$, and $H$ in Section~\ref{sec:energy} are readily discretized,
as are the force densities in~\ref{sec:forces}. Application of the dense discrete differential operators is
performed in parallel with a parallel implementation of \texttt{BLAS}. Lagrangian forces can therefore be computed
in parallel with few thread synchronizations.

We now have a method for computing a suitable set of points and for discretizing $\F$ for use in~%
\eqref{eq:ib-spread}. To compute a force from a force density, we need to compute the set of quadrature weights
$\weight[j]$, $j=1, \ldots, \sample\cardinality$. The following section is devoted to describing our method for
computing these weights.

\subsection{RBF-based quadrature}\label{sec:rbf-quadrature}

Our approach to generating quadrature weights $\weight[j]$ for numerical integration on an immersed structure
$\interface$ (such as RBCs and platelets) at the spreading points $\Xsj$ is to first precompute a fixed set of
quadrature weights $\omega_j$ on $\sphere$, then map these weights to $\interface$ during each time-step. For the
following discussion, note that $\weight[j](\theta, \varphi, t)$ is technically a function of time $t$, since
$\interface$ moves in time. However, for brevity, we suppress the arguments for the weights.

First, we discuss how to compute RBF-based meshless quadrature weights on $\sphere$. Given a function
$\psi:\sphere\to\R$ and a radial kernel $\phi:\sphere\times\sphere\to\R$, we now wish to find a set of quadrature
weights $\omega_j$ such that
\begin{equation}\label{eq:quad-desire}
    \int\limits_{\sphere} \psi(\Xp) \d\Xp \approx \sum_{j=1}^{\sample\cardinality} \omega_j \psi(\Xpsj).
\end{equation}
We use a variant of the technique described by Fuselier \latin{et al.}~\cite{Fuselier:2013coba}.  Choosing
$\psi(\Xp) = \phi(\|\Xp-\sample\Xp_i\|)$ for each $\sample\Xp_i$,~\eqref{eq:quad-desire} becomes
\begin{equation}
    \sum_{j=1}^{\sample\cardinality} \omega_j \phi(\|\Xpsj-\sample\Xp_i\|)
    \approx \int\limits_{\sphere}\phi(\|\Xp-\sample\Xp_i\|) \d\Xp := \L\phi|_{\sample\Xp_i}.
    \label{eq:cond1t}
\end{equation}
However, because
\begin{equation}
    \|\Xp-\sample\Xp_i\|^2
        = (\Xp-\sample\Xp_i)\cdot(\Xp-\sample\Xp_i)
        = 2 - 2\Xp\cdot\sample\Xp_i
\end{equation}
depends only on the angle between the two points, $\L\phi$ is constant over the sphere, and we denote this constant $-I_\phi$. We require further that the $\omega_j$ sum to the
surface area of $\sphere$, \latin{i.e.},
\begin{equation}\label{eq:weight-sum}
    \sum_{j=1}^{\sample\cardinality} \omega_j  = 4\pi.
\end{equation}
Treating $I_\phi$ as an unknown scalar, we rewrite the constraints~\eqref{eq:cond1t} and~\eqref{eq:weight-sum} in
the symmetric block linear system
\begin{equation}\label{eq:rbf-quadrature}
    \begin{bmatrix}
            \Phi & \arr{1} \\ \arr{1}^T & 0
    \end{bmatrix}\begin{bmatrix}
            \arr{\omega} \\ I_\phi
    \end{bmatrix} = \begin{bmatrix}
            \arr{0} \\ 4\pi
    \end{bmatrix},
\end{equation}
where $\arr{\omega}$ are the unknown weights, $\Phi$ represents the evaluations of $\phi$, and $\arr{0}$ and
$\arr{1}$ are vectors of $\sample\cardinality$ zeros and ones, respectively. $I_{\phi}$ serves as a Lagrange
multiplier that enforces~\eqref{eq:weight-sum}, and $-I_\phi$ is a good approximation to $\L\phi$. By choosing
$\phi(r) = r$, we guarantee a unique solution and obtain weights that when used in a quadrature rule, converge at 3\textsuperscript{rd} order in the sample site spacing. It is possible to improve the order of the quadrature weights by increasing the order
of the polyharmonic spline RBF at the potential cost of poorer conditioning and either loss of invertibility or requiring
knowledge of higher-order moments~\cite{Fuselier:2013coba}.

Once the quadrature weights $\omega_j$, $j=1, \ldots, \sample\cardinality$ are computed, it is straightforward to
compute the corresponding quadrature weights $\weight[j]$, $j=1, \ldots, \sample\cardinality$ for integration at
the spreading points $\sample\X \in \interface(t)$. Recall that $\interface$ has parametrization
$\X(\theta, \varphi, t)$ and Jacobian $J(\theta, \varphi)$ for a given time t. This implies that the point
$\X(\theta, \varphi)$ directly corresponds to to $\Xp(\theta, \varphi)$. Using this fact and the fact that the
determinant of the Jacobian for the spherical coordinate mapping (for radius 1) is $\cos\varphi$, we can use a
change of variables to express the infinitesimal area $\d\X$ on $\interface$ as
\begin{equation}\label{eq:quad-cov}
    \d\X
    = J(\theta, \varphi)\d\qs
    = J(\theta, \varphi)\sec\varphi\d\Xp,
\end{equation}
where $\d\qs$ is an infinitesimal area in parameter space. Note that the weights $\omega_j$ above are discrete
analogs of $\d\Xp$ at $\Xpsj$. The discrete analog of $\d\qs$ at the $j^\text{th}$ sample site,
\begin{equation}
    \sigma_j=\sec\varphi_j\omega_j,
\end{equation}
can therefore be computed at the outset of a simulation. To avoid numerical issues, we require that
$\cos\varphi_j \neq 0$ for each sample site. This is true everywhere on $\sphere$ except the poles,
$(0, 0, \pm1)$, which the Bauer spiral~\eqref{eq:bauer-spiral} conveniently avoids. Finally, the quadrature weight
$\weight[j]$ for $\interface$ at the spreading point $\Xsj$ is given by
\begin{equation}
    \weight[j] = \sigma_j J_j,
\end{equation}
where $J_j = J(\theta_j, \varphi_j)$. Computing $\weight[j]$ given $\sigma_j$ and $J_j$ amounts to a single
multiplication, which can be done trivially in parallel. This produces a set of (time-varying) weights
$\weight[j]$, $j=1, \ldots, \sample\cardinality$ that allows us to integrate functions on $\interface(t)$ when
their samples are given at $\Xsj$, $j=1, \ldots, \sample\cardinality$, as required by~\eqref{eq:ib-spread-disc}.

\section{Results}\label{sec:results}

We have taken a number of departures from the traditional application of the IB method and RBC models. In the
following sections, we establish convergence for our implementation of RBF-IB, demonstrate its energy dissipation properties, demonstrate impermeability of the
RBCs, and observe typical RBC behaviors before presenting the results of our whole blood simulations.

\subsection{Convergence study}\label{sec:convergence}

In this section, we perform a series of tests on a single perturbed RBC undergoing relaxation. We simplify the RBC
model and use only Skalak's Law. We expect the IB method to approximate the fluid velocity at first order for thin
shells, as it cannot recover the pressure jump across the interface. We stretch the RBC by a factor of 1.1 in the
$z$ direction and compress it in the $x$ direction to maintain its reference volume. We place the cell in the
center of a $16\um\times16\um\times16\um$ domain with homogeneous Dirichlet boundary conditions in the $y$
direction and periodic boundaries elsewhere. The fluid velocity is initially zero. The cell is then allowed to
relax for $180\us$.

We use the 3-point kernel $\hat{\Dirac}_1$ derived by Roma \latin{et al.}~\cite{Roma:1999tx} for spreading and
interpolation and the 2-stage RK method described in \cref{sec:ib} to advance the fluid velocity. Fluid grids are
chosen to have $20r$ grid points per $16\um$ in each direction for $r$ from 1 to 6. On successive grids, we
compare the fluid velocity at cell centers in a regular grid of $20^3$ cells and surface positions at 1000 surface
points. For convergence of order $p$, we expect the ratio of successive differences to satisfy
\begin{equation}
    \frac{\epsilon_r}{\epsilon_{r+1}} = \left|\frac{{((r+1)/r)}^p-1}{{((r+2)/(r+1))}^{-p}-1}\right|,
\end{equation}
which we solve numerically to determine $p$. We compute the differences in $\X$ using discrete versions of the $L^2$
and $L^\infty$ norms,
\begin{gather}
    \|\X(\theta, \varphi)\|_2^2 =
    \int\limits_{\sphere} \X(\theta, \varphi)\cdot\X(\theta, \varphi) \d\qs \quad\text{and} \\
    \|\X(\theta, \varphi)\|_\infty^2 =
    \max_{(\theta, \varphi)} \X(\theta, \varphi)\cdot\X(\theta, \varphi),
\end{gather}
respectively. In all cases, we label these differences as ``errors'', though this is not strictly true.

\Cref{tab:u-rbc-conv,tab:x-rbc-conv} show $L^2$ and $L^\infty$ errors in $\u$ and $\X$ between successive grids. We observe first-order convergence, as expected. 

\begin{table}[tb]
    \centering
    \caption[Convergence of fluid velocities for relaxing RBC test]{%
Convergence of $\u$ for a sequence of grids. The refinement ratio $r$, defined as the refinement in $h$ relative
to the coarsest grid, determines the simulation parameters: $rh = 0.8\um$ and $r\timestep = 180\ns$. Differences are
computed between grids of refinement factor $r$ and $r+1$. Values of $\u$ are sampled at $t = 180\us$ at cell
centers on the coarsest grid.
    }\label{tab:u-rbc-conv}
    \begingroup
    \setlength{\tabcolsep}{9pt}
    \renewcommand{\arraystretch}{1.5}
    \begin{tabular}{c|cc|cc}
                                                                                     \toprule
        $r$ & $L^2$ error            & order   & $L^\infty$ error       & order   \\ \midrule
        1   & $\scinot{1.74592}{-3}$ &         & $\scinot{1.59995}{-2}$ &         \\
        2   & $\scinot{9.92788}{-5}$ &         & $\scinot{6.52038}{-4}$ &         \\
        3   & $\scinot{3.65264}{-5}$ & 6.85983 & $\scinot{3.34322}{-4}$ & 1.06075 \\
        4   & $\scinot{2.31069}{-5}$ & 2.87612 & $\scinot{2.30732}{-4}$ & 1.60689 \\
        5   & $\scinot{1.65898}{-5}$ & 1.25668 & $\scinot{2.28193}{-4}$ & 0.83030 \\ \bottomrule
    \end{tabular}
    \endgroup
\end{table}

\begin{table}[tb]
    \centering
    \caption[Convergence of surface positions for relaxing RBC test]{%
Convergence of $\X$ for a sequence of grids. The refinement ratio $r$, defined as the refinement in $h$ relative
to the coarsest grid, determines the simulation parameters: $r\timestep = 180\ns$, $n_d = 125r^2$, and
$n_s = 500r^2$. Differences (errors) are computed between grids of refinement factor $r$ and $r+1$. Values of $\X$ are sampled
at $t = 180\us$ at $1000$ Bauer spiral points.
    }\label{tab:x-rbc-conv}
    \begingroup
    \setlength{\tabcolsep}{9pt}
    \renewcommand{\arraystretch}{1.5}
    \begin{tabular}{c|cc|cc}
                                                                                     \toprule
        $r$ & $L^2$ error            & order   & $L^\infty$ error       & order   \\ \midrule
        1   & $\scinot{6.05447}{-3}$ &         & $\scinot{2.68249}{-3}$ &         \\
        2   & $\scinot{1.61678}{-3}$ &         & $\scinot{6.86140}{-4}$ &         \\
        3   & $\scinot{7.69150}{-4}$ & 2.74664 & $\scinot{3.30734}{-4}$ & 2.86457 \\
        4   & $\scinot{4.66203}{-4}$ & 1.89474 & $\scinot{1.86139}{-4}$ & 1.84435 \\
        5   & $\scinot{2.82578}{-4}$ & 1.46574 & $\scinot{1.18713}{-4}$ & 1.82742 \\ \bottomrule
    \end{tabular}
    \endgroup
\end{table}

\subsection{Energy Estimates}\label{sec:energy-est}

As mentioned in \cref{sec:rbfib}, the RBF-IB method uses different sets of Lagrangian points for spreading forces
and updating immersed structures. We now track the total energy of the combined fluid-RBC system in the relaxation
test just described.  This energy is given as 
\begin{align}
    E = \int\limits_\domain \left[\frac{\density}2\u\cdot\u + \int_\interface \Dirac(\x-\X)W_\text{Sk}(\X, \ldots)\d\X\right] \d\x.
\label{eq:energy}
\end{align}
The first term corresponds to the energy of the fluid, and the second to the energy of the structure. Since the
relaxation test involves an initial increase in kinetic energy of the fluid followed by a gradual dissipation of
the potential energy of the RBC as it relaxes according to the Skalak law, we expect the total energy of the
system to decrease over time. We plot a discretized version of the energy~\eqref{eq:energy} as a function of time
in \cref{fig:energies} for refinement factors $r=5$ and $r=6$. \Cref{fig:energies} clearly shows that our method
dissipates energy in the manner expected in this problem. We also observed that the RBF-IB method was stable for
representing RBCs and platelets in whole blood simulations.

\begin{figure}[tb]
\centering
\begin{tikzpicture}
\begin{groupplot}[
    group style={
        y descriptions at=edge left,
        group name=energy,
        group size=2 by 1
    },
    width=2.5in,
    height=2.5in,
    xmin=-10,
    xmax=190,
    ymin=9e-12,
    ymax=1.1e-9,
    ymode=log,
    log basis y=10,
    log origin=infty,
    axis x line=bottom,
    axis y line=left,
    xlabel={time ($\us$)},
    xlabel near ticks,
    ylabel near ticks,
]
\nextgroupplot[ylabel={energy ($\erg$)}]
    \addplot+[only marks, mark options={fill=tol/vibrant/blue}] coordinates {%
        (  0, 48.4709e-11)
        ( 18, 7.56723e-11)
        ( 36, 3.71642e-11)
        ( 54, 3.14419e-11)
        ( 72, 3.0141e-11)
        ( 90, 2.94874e-11)
        (108, 2.89456e-11)
        (126, 2.84339e-11)
        (144, 2.79377e-11)
        (162, 2.74534e-11)
        (180, 2.69799e-11)
    }; \label{plot:energy100}
    \node [fill=white] at (rel axis cs: 0.075, 0.95) {(a)};
\nextgroupplot
    \addplot+[only marks, mark options={fill=tol/vibrant/blue}] coordinates {%
        (  0, 4.84709e-10)
        ( 18, 7.30549e-11)
        ( 36, 3.64163e-11)
        ( 54, 3.122e-11)
        ( 72, 3.00294e-11)
        ( 90, 2.93954e-11)
        (108, 2.88559e-11)
        (126, 2.83433e-11)
        (144, 2.78457e-11)
        (162, 2.736e-11)
        (180, 2.68852e-11)
    }; \label{plot:energy120}
    \node [fill=white] at (rel axis cs: 0.075, 0.95) {(b)};
\end{groupplot}
\end{tikzpicture}
\caption{%
Energy~\eqref{eq:energy} as a function of time for the relaxing RBC test of \cref{sec:convergence} with refinement
factors (a) $r=5$ and (b) $r=6$. The refinement factor determines the simulation parameters: spacestep
$rh = 0.8\um$, timestep $r\timestep=180\ns$, $\data\cardinality=125r^2$, and $\sample\cardinality=500r^2$.
}\label{fig:energies}
\end{figure}

\subsection{RBC tumbling and tank-treading}

\begin{figure}[tb]
    \centering
    \begin{subfigure}{\textwidth}
    \begin{minipage}{0.2\textwidth}
        \centering
        \topinset{(a)}{\includegraphics[trim=75 100 75 100, clip, width=\textwidth]{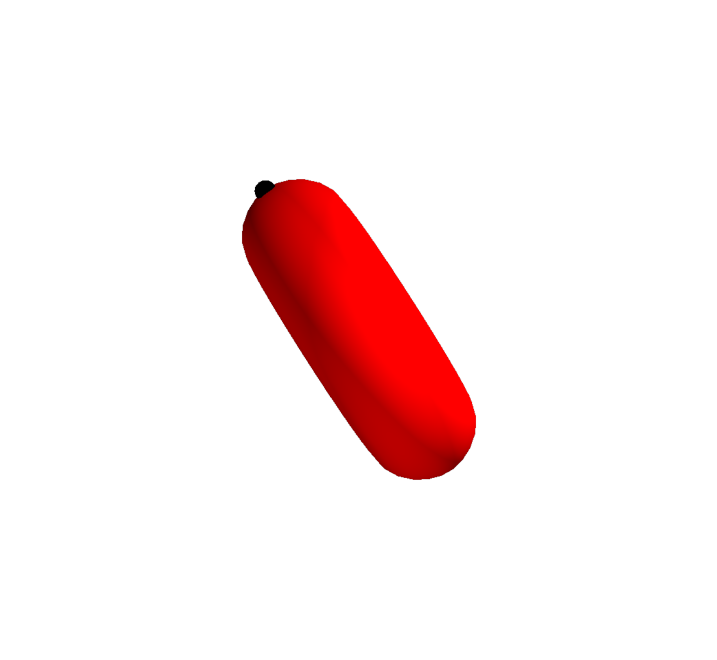}}{0.5cm}{0.25cm}\\
        $\dot{\gamma}t = 0$
    \end{minipage}%
    \begin{minipage}{0.2\textwidth}
        \centering
        \includegraphics[trim=75 100 75 100, clip, width=\textwidth]{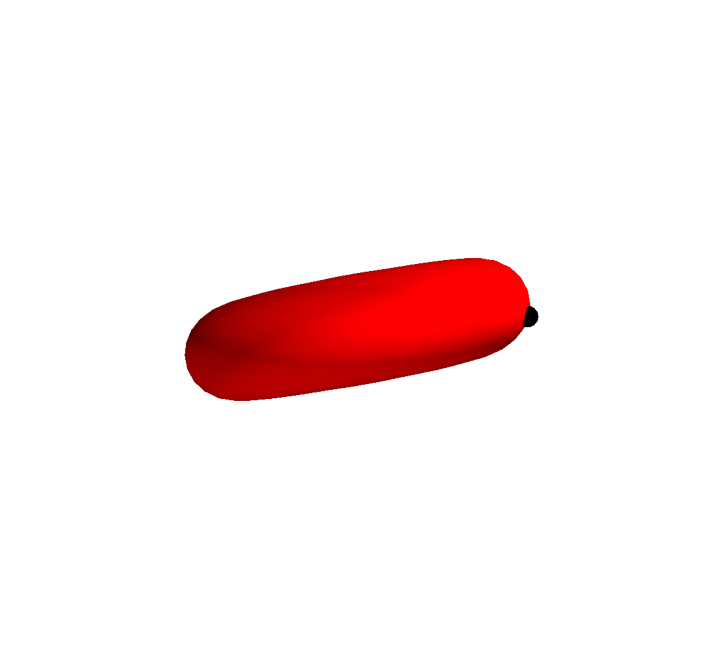}\\
        $\dot{\gamma}t = 5$
    \end{minipage}%
    \begin{minipage}{0.2\textwidth}
        \centering
        \includegraphics[trim=75 100 75 100, clip, width=\textwidth]{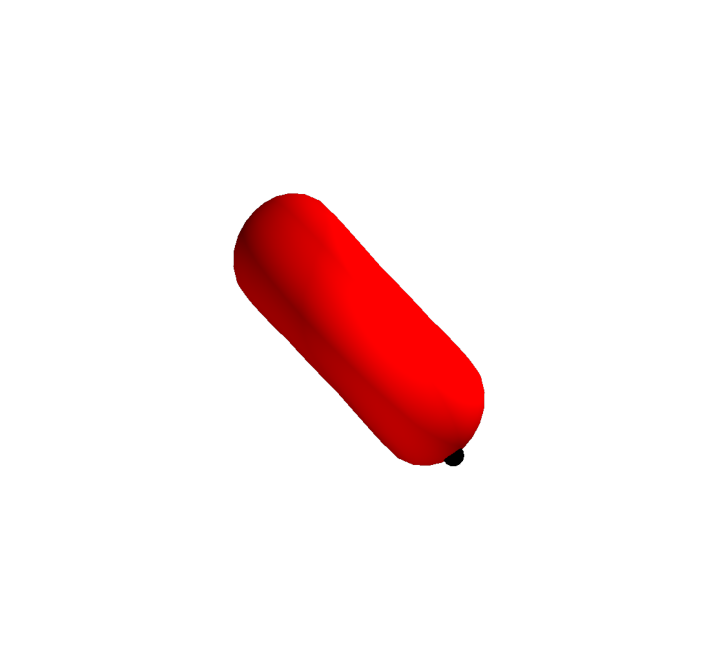}\\
        $\dot{\gamma}t = 10$
    \end{minipage}%
    \begin{minipage}{0.2\textwidth}
        \centering
        \includegraphics[trim=75 100 75 100, clip, width=\textwidth]{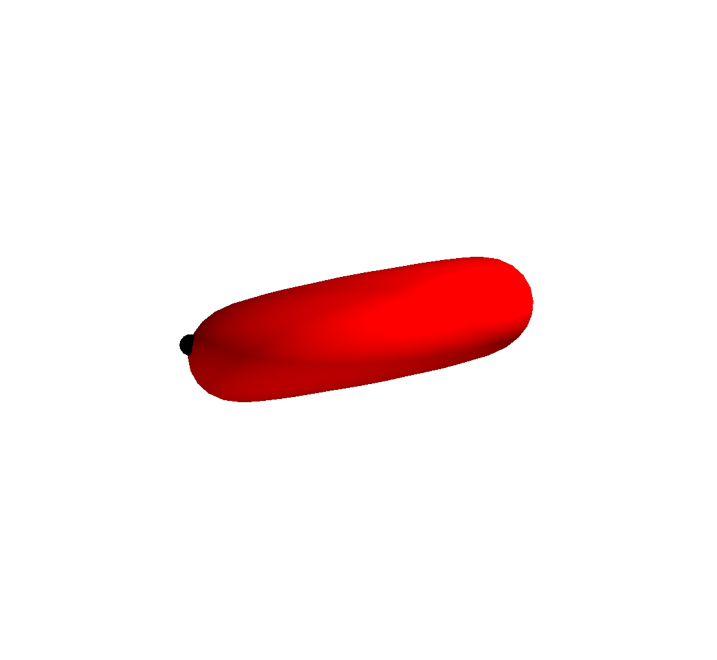}\\
        $\dot{\gamma}t = 15$
    \end{minipage}%
    \begin{minipage}{0.2\textwidth}
        \centering
        \includegraphics[trim=75 100 75 100, clip, width=\textwidth]{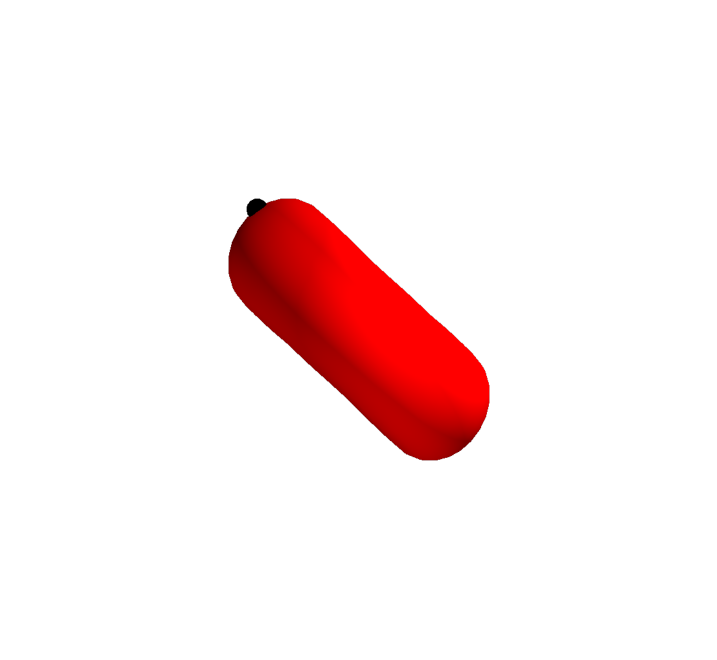}\\
        $\dot{\gamma}t = 20$
    \end{minipage}%
    \phantomsubcaption\label{fig:tumble}
    \end{subfigure}
    \begin{subfigure}{\textwidth}
    \begin{minipage}{0.2\textwidth}
        \centering
        \topinset{(b)}{\includegraphics[trim=75 125 75 75, clip, width=\textwidth]{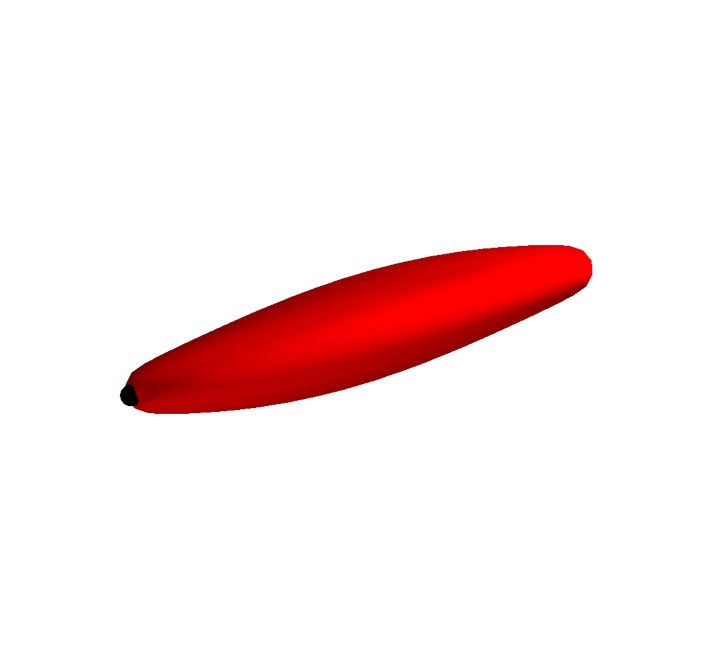}}{0.5cm}{0.25cm}\\
        $\dot{\gamma}t = 19$
    \end{minipage}%
    \begin{minipage}{0.2\textwidth}
        \centering
        \includegraphics[trim=75 125 75 75, clip, width=\textwidth]{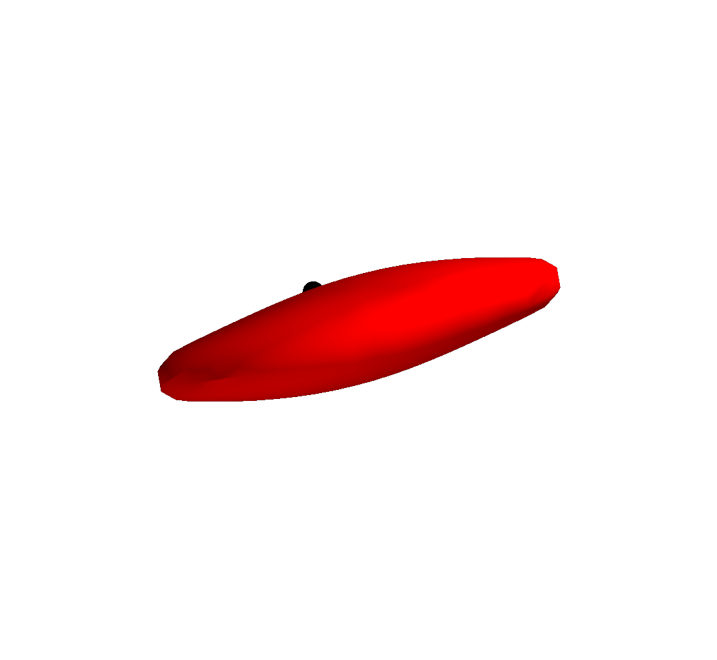}\\
        $\dot{\gamma}t = 26$
    \end{minipage}%
    \begin{minipage}{0.2\textwidth}
        \centering
        \includegraphics[trim=75 125 75 75, clip, width=\textwidth]{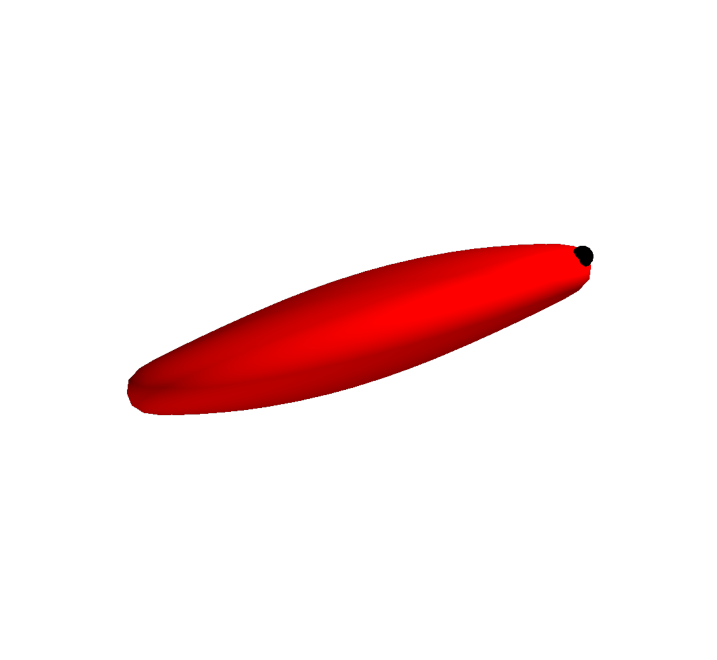}\\
        $\dot{\gamma}t = 33$
    \end{minipage}%
    \begin{minipage}{0.2\textwidth}
        \centering
        \includegraphics[trim=75 125 75 75, clip, width=\textwidth]{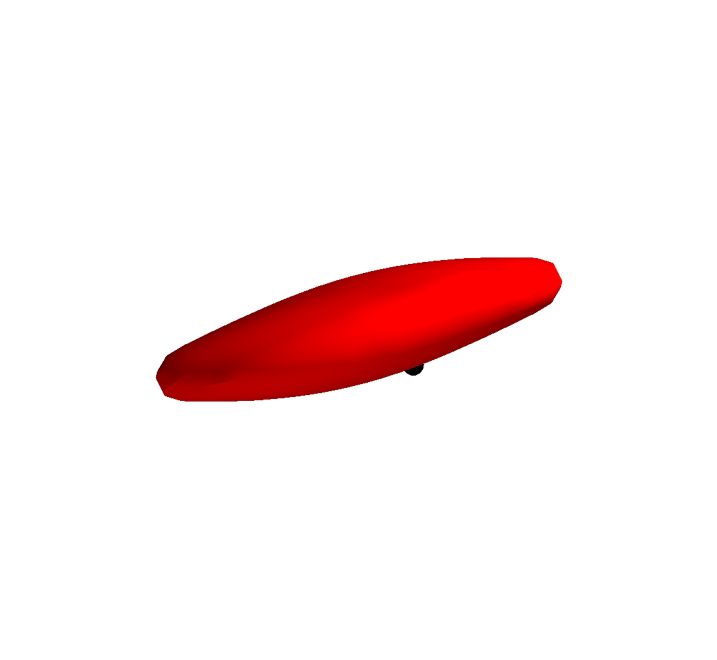}\\
        $\dot{\gamma}t = 40$
    \end{minipage}%
    \begin{minipage}{0.2\textwidth}
        \centering
        \includegraphics[trim=75 125 75 75, clip, width=\textwidth]{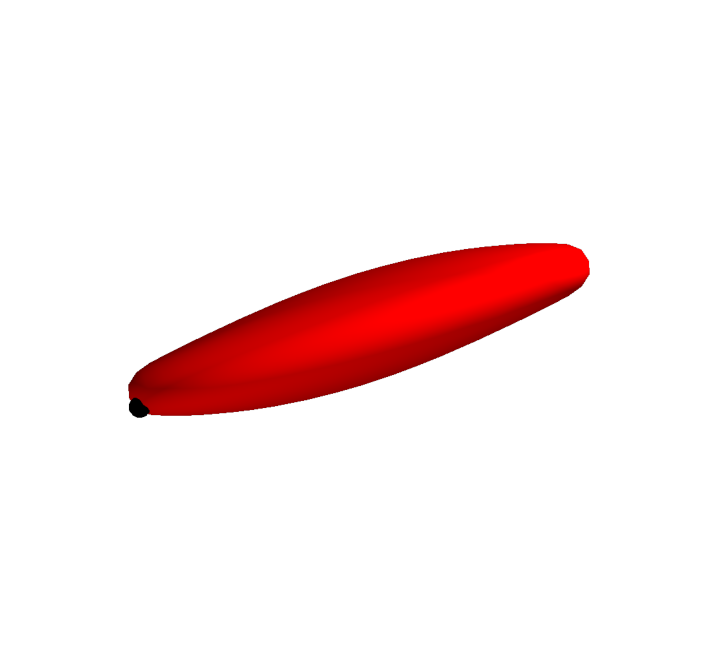}\\
        $\dot{\gamma}t = 47$
    \end{minipage}%
    \phantomsubcaption\label{fig:tread}
    \end{subfigure}
    \caption{%
Our model RBC exhibits (a) a tumbling behavior under low shear ($\dot{\gamma} = 50\si{\per\second}$) conditions
and (b) tank-treading under high shear ($\dot{\gamma} = 1000\si{\per\second}$) conditions.
    }%
    \label{fig:tumble-tread}
\end{figure}

Few RBC models include viscoelastic forces. Fedosov \latin{et al.}~\cite{Fedosov:2010bc} use a particle-based
method to simulate the fluid and cells, where particles representing the RBC membrane experience drag and random
forces from relative motion between RBC particles. Gounley \& Peng~\cite{Gounley:2015ho} use the IB machinery to
spread membrane viscosity to the fluid, thereby modifying the fluid stress term. Our intent in adding a
viscoelastic response is to aid in the numerical stability of the discrete RBCs. We wish to verify that the
extended model, with dissipative force, retains the ability to tumble and tank-tread, which is commonly used to
validate other RBC models~\cite{Yazdani:2011cl, Omori:2012hw, Fai:2013do, Xu:2013kk}. To that end, we place a
single RBC with $\data\cardinality=625$ and $\sample\cardinality = 2500$ in the same domain as the previous
section, now discretized to have $h = 0.4\um$ and with moving top and bottom walls. In the interest of reducing
simulation time, we now use the backward-forward Euler timestepping scheme with a time step of $\timestep=0.1\us$.
Here, the IB interaction operations use the 4-point B-spline, $\kernel(r) = B_3(r)$, which was first considered by
Lee~\cite{Lee:2020tf}. It is similar in shape to the Roma kernel but has better smoothness properties. To recover
the tumbling motions, the top wall has a fixed velocity of $\u_b=400\umpersec$ and the bottom wall $-\u_b$. This
generates a shear rate of $\shear=50\persec$ in the absence of cells. For tank-treading experiments, we use
$\u_b=8\mmpersec$ to generate a shear rate of $\shear=1000\persec$. These values are chosen outside the
transitional region between tumbling and tank-treading for the elastic parameters used for the RBC~%
\cite{Kruger:2013ji}. The velocity field is initially steady for flow without cells. We rotate the cell 1 radian
about the $x$-axis from a horizontally aligned orientation and place it at the center of the domain. The RBC
exhibits tumbling for $\shear=50\persec$ and tank-treading for $\shear=1000\persec$. \Cref{fig:tumble-tread} shows
one period of each.

\subsection{Collision tests}

With whole blood simulation as the ultimate goal, we must ensure that the method can effectively capture cell-cell
interactions. The RBF-IB method has been applied to flow around multiple platelets in an aggregate in 2D simulations, but those
cells were kept apart by a network of springs~\cite{Shankar:2015km}. To study the interaction between cells, we
devise a series of tests in which we force two RBCs to collide. The aim is to verify that the cells remain
distinct. Using too few data sites could allow the cells to come too close to one another or even to interpenetrate, as the regularized delta function $\Dirac_h$ then causes them to be treated as a single unit. Cells that ``fuse'' in this manner are problematic,
generally causing the simulation to end (due to stability-related time-step restrictions) when the cells attempt to separate.

\begin{figure}[tb]
    \centering
    \begin{subfigure}[t]{.25\textwidth}
        \centering
        \topinset{(a)}{\includegraphics[width=\textwidth]{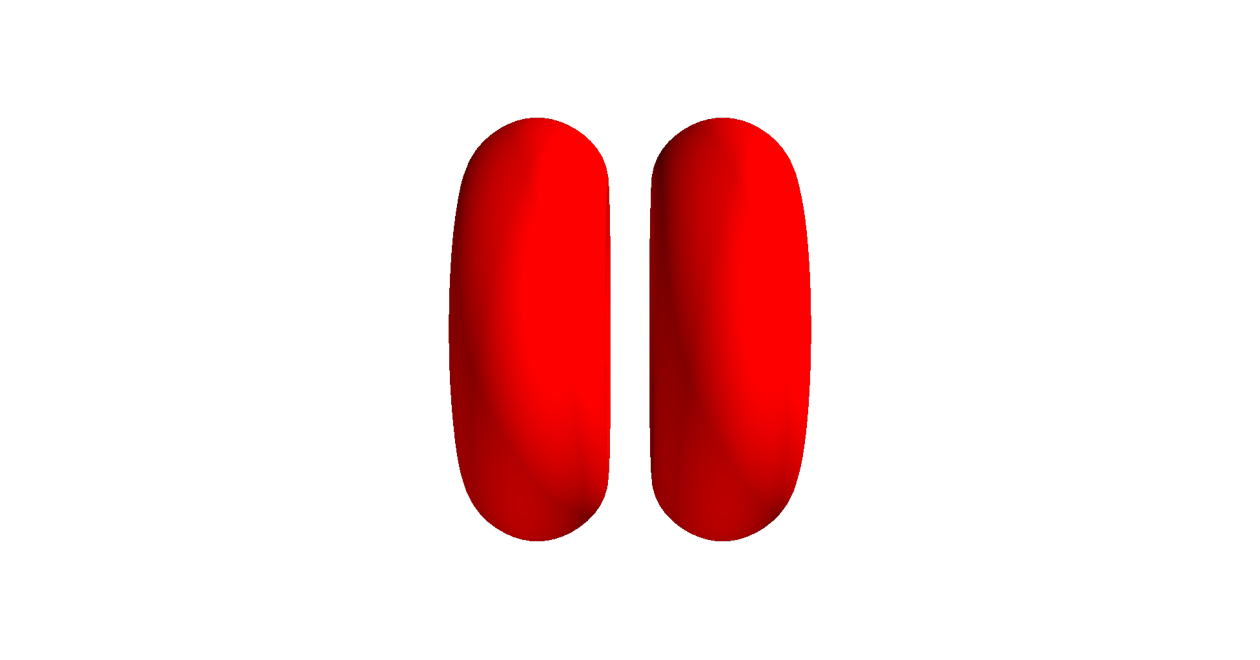}}{0.125cm}{0.25cm} \\
        $t = 0\ms$
    \end{subfigure}%
    \begin{subfigure}[t]{.25\textwidth}
        \centering
        \topinset{(b)}{\includegraphics[width=\textwidth]{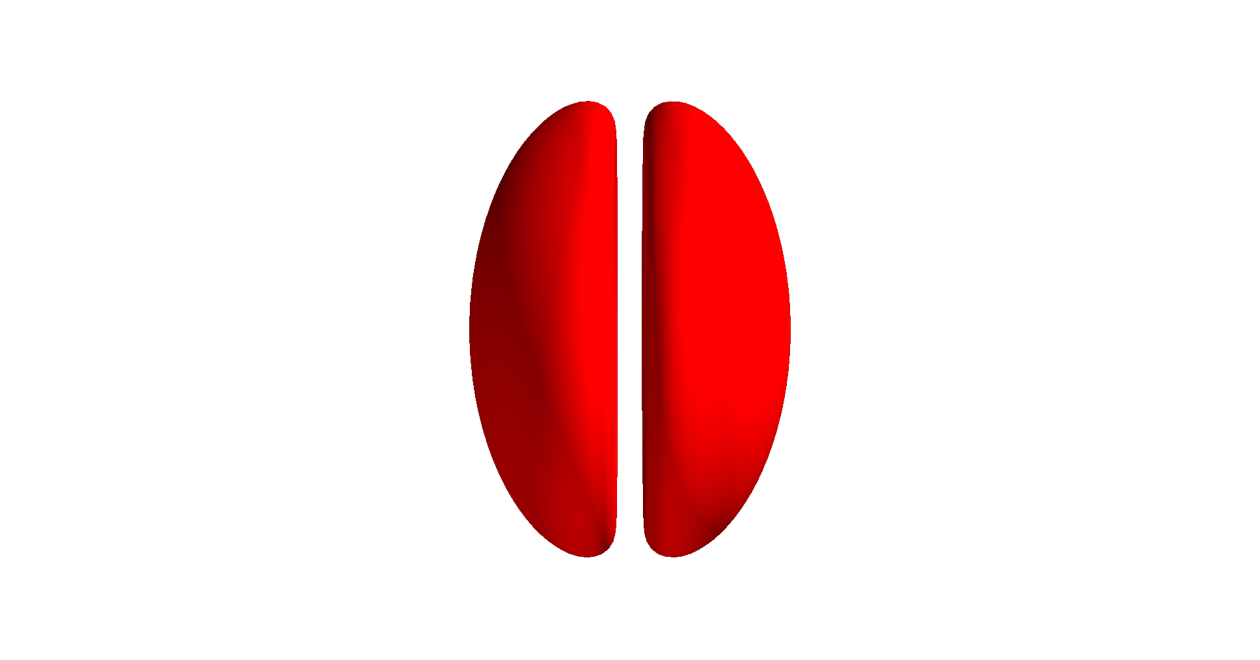}}{0.125cm}{0.25cm} \\
        $t = 1.5\ms$
    \end{subfigure}%
    \begin{subfigure}[t]{.25\textwidth}
        \centering
        \topinset{(c)}{\includegraphics[width=\textwidth]{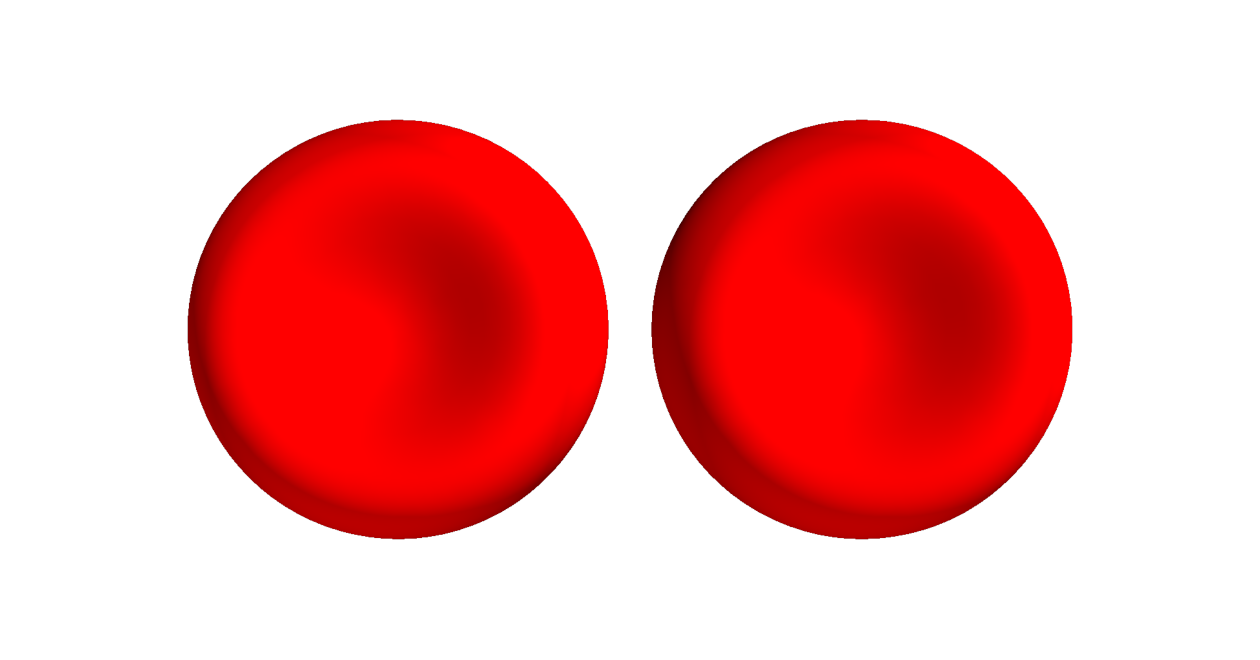}}{0.125cm}{0.25cm} \\
        $t = 0\ms$
    \end{subfigure}%
    \begin{subfigure}[t]{.25\textwidth}
        \centering
        \topinset{(d)}{\includegraphics[width=\textwidth]{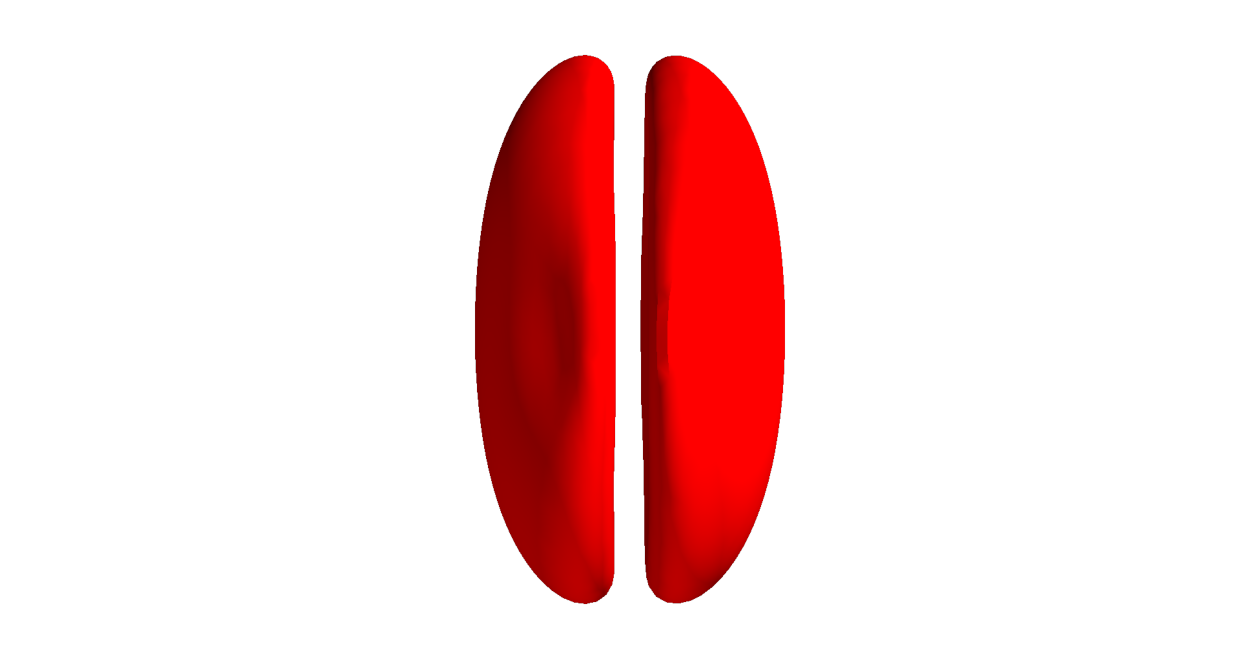}}{0.125cm}{0.25cm} \\
        $t = 1.1\ms$
    \end{subfigure}

    \vspace{1em}

    \begin{subfigure}[t]{.25\textwidth}
        \centering
        \topinset{(e)}{\includegraphics[width=\textwidth]{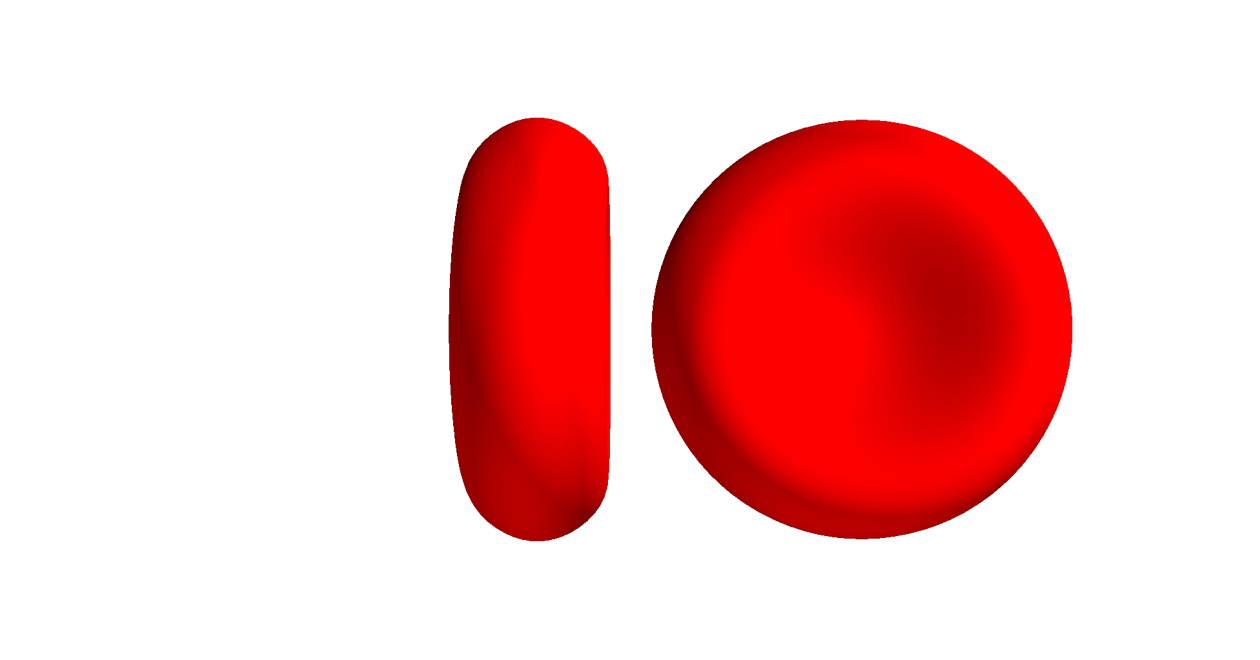}}{0.125cm}{0.25cm} \\
        $t = 0\ms$
    \end{subfigure}%
    \begin{subfigure}[t]{.25\textwidth}
        \centering
        \topinset{(f)}{\includegraphics[width=\textwidth]{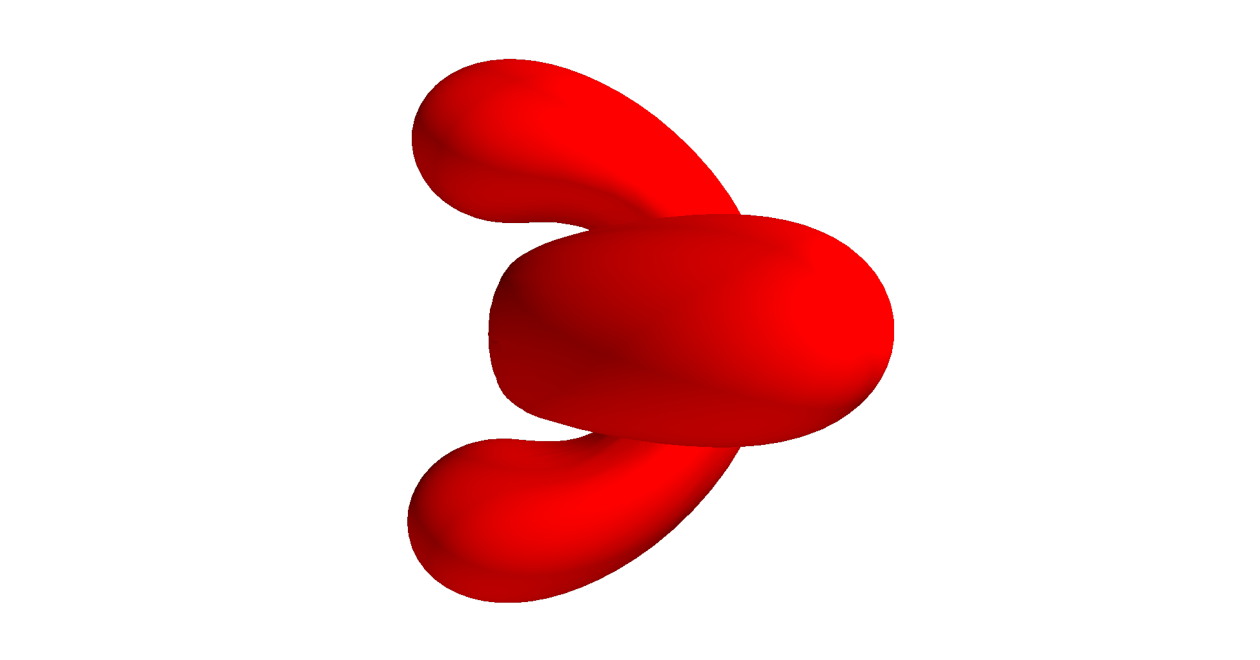}}{0.125cm}{0.25cm} \\
        $t = 1.5\ms$
    \end{subfigure}%
    \begin{subfigure}[t]{.25\textwidth}
        \centering
        \topinset{(g)}{\includegraphics[width=\textwidth]{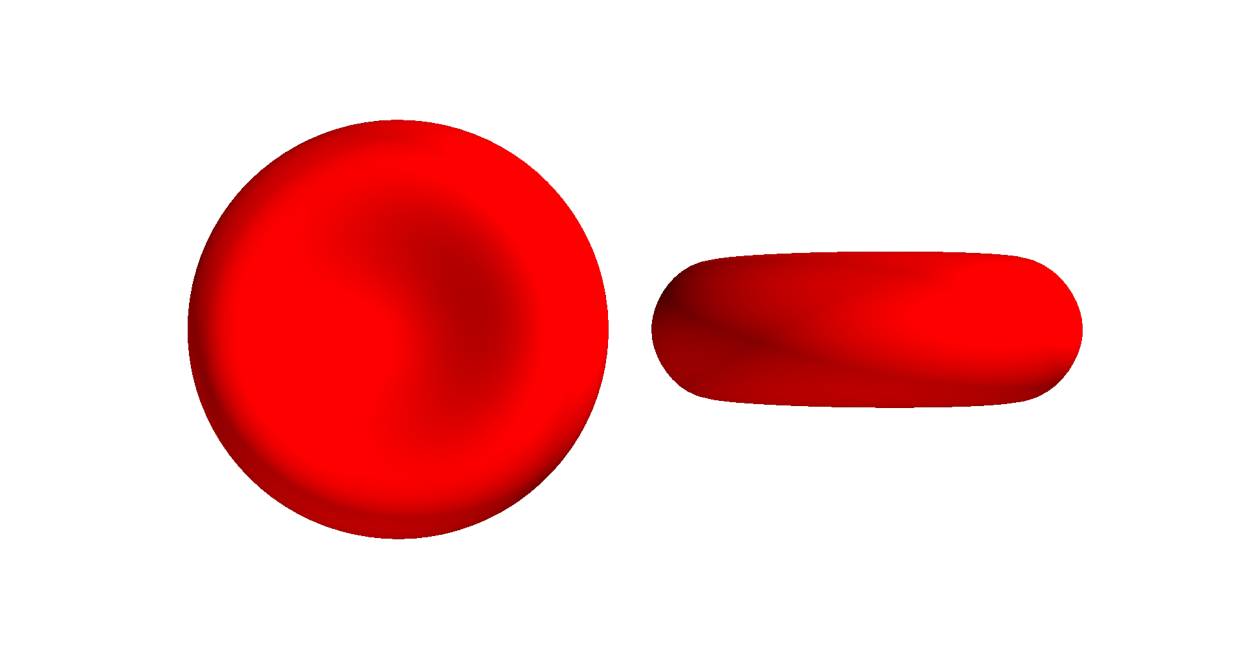}}{0.125cm}{0.25cm} \\
        $t = 0\ms$
    \end{subfigure}%
    \begin{subfigure}[t]{.25\textwidth}
        \centering
        \topinset{(h)}{\includegraphics[width=\textwidth]{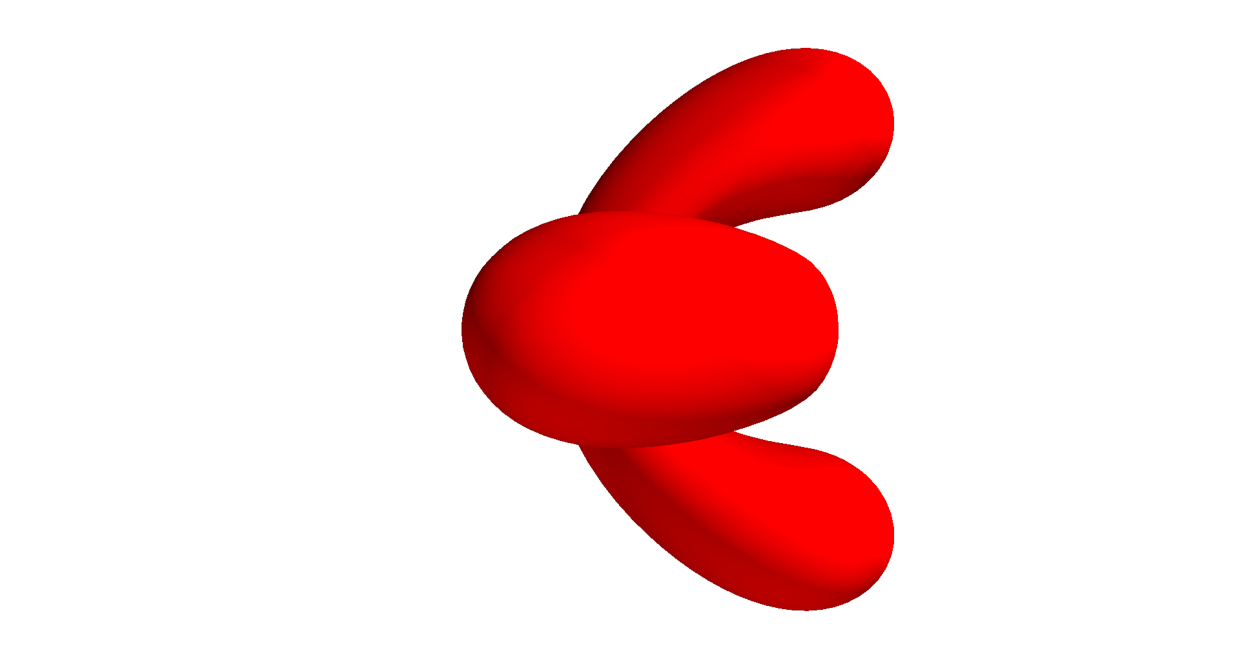}}{0.125cm}{0.25cm} \\
        $t = 1.1\ms$
    \end{subfigure}
    \caption{%
Collision tests between two RBCs. A fictitious force is added to the RBCs to draw them together. (a) and (b) The
RBCs are initially aligned with concavities facing one another. By $1.5\ms$, the cells take on a hemispherical
shape. The concavities in the gap are maintained. Shortly thereafter, asymmetries in the setup lead to the cells
sliding past one other. (c) and (d) The RBCs are initially aligned with their edges facing one another. By
$1.1\ms$, the cells take on a hemispherical shape. The remnants of the concavity can be seen on the left cell in
(d). Shortly thereafter, these cells also slide past one another. (e) and (f) The RBCs are initially aligned with
the edge of one facing a concavity of the other. The cells wrap around each other by $1.5\ms$, taking on a bulbous
banana shape. (g) and (h) The RBCs are initially aligned with their edges facing one another with one of the cells
rotated about the axis $\e_1+\e_3$ by $\pi/2$. By $1.1\ms$, the cells wrap around each other, again taking on the
bulbous banana shape.
    }%
    \label{fig:collisions}
\end{figure}

We use the same physical domain as in the previous two sections, now with $h = 0.2\um$, and place two
RBCs therein, each with $\data\cardinality=2500$ and $\sample\cardinality=10000$. The ratio
$\sample\cardinality/\data\cardinality=4$ is chosen so that sample sites with spacing $h$ means data sites have
spacing $2h$. We believe this to suffice in preventing cells from intersecting, but this is not guaranteed if the
points do not maintain appropriate spacing throughout a simulation. We use the 2-stage RK method with time step
$\timestep = 50\ns$. To interpolate velocities and spread forces, we use the 4-point cosine $\kernel$~%
\cite{Peskin:2002go}. The cells are placed with cell centers on the line $x = z$, $y = 8\um$. They are initially
separated by a gap of $4h = 0.8\um$ between their convex hulls, \latin{i.e.}, ignoring the concavities. Inspired
by Crowl \& Fogelson~\cite{Erickson:2011cf}, we add the fictitious force density
\begin{equation}
    \F_\text{fict} = \pm 0.1\dynpercm\cdot(\e_1+\e_3)/\sqrt{2},
\end{equation}
to each cell, where the sign is chosen so the force points into the gap, to draw the cells together. Success in
these tests implies that this configuration of data and sample sites, the spatial resolution, and the time step
are acceptable for whole blood simulations.

Initial conditions and configurations after a short time are illustrated in \cref{fig:collisions}, where we view
them from above the $x=z$ plane. In each case, the cells move slightly closer together and then undergo
considerable deformation. The data sites are initially approximately $2h$ apart from each other. No problems seem
to arise from this, and in some cases the cells eventually attempt to slide past one another.  We also deduce that
the IB method with the cosine kernel can resolve interactions at a distance of $h$ to $2h$. We consider cells
passing within this threshold to be in contact. Throughout the simulation, the cells remain distinct, and the
simulations end due to extreme forces triggering the stopping condition~\cite{Agresar:1998wv}
\begin{equation}
    \timestep > \frac14\sqrt{\frac{h\rho}{\|\f\|_\infty}}.
\end{equation}
For the remainder of this section, we consider only this arrangement of data and sample sites and this grid
resolution.

\subsection{Whole blood}\label{sec:whole-blood}

In what follows, we consider a $16\um\times12\um\times16\um$ domain with periodic boundaries in the $x$ and $z$
directions and with Dirichlet boundary conditions in the $y$ direction. The fluid velocity is initially zero
except at the top boundary, where it moves at $12\mmpersec$. In the absence of cells, the flow tends toward steady
Couette flow with a shear rate of $\shear=1000\persec$. This serves as our model near-wall region of a blood
vessel.

For whole blood simulations, we return to the 4-point B-spline, $B_3$, as the IB kernel.  Because some RBC
configurations generate large forces, the stopping condition of the previous section limits us to small timesteps.
In the interest of reducing simulation time, we use the backward-forward Euler scheme with $\timestep=50\ns$. In
tests, we observe qualitative agreement between this scheme and the two-stage RK scheme. We have already settled
on an RBC discretization in the previous section. We use the same spiral method to discretize the platelet, but
with 900 data and sample sites. Using the same number of sample sites and data sites aligns more closely with
traditional IB methods. We also find that the Bauer spiral places points more densely along the edge of the
platelet, which is helpful in resolving the large curvatures there. We parametrize the surface of the endothelium
over $(\theta, \varphi)\in{[0, 2\pi)}^2$ with reference shape
\begin{equation}
    \vec{\hat{X}}_\text{endo}(\theta, \varphi) = \begin{bmatrix}
            16\um\cdot(\theta/2\pi)  \\
            y(\theta, \varphi) \\
            16\um\cdot(\varphi/2\pi)
    \end{bmatrix},
\end{equation}
where $y(\theta, \varphi)$ depends on the shape under study. The endothelium is discretized using 16000 points,
defined by the spiral
\begin{equation}
    \begin{aligned}
    \varphi_i &= 2\pi (i-1)/N, \\
    \theta_i &= \modulo\left(\left\lceil\!\sqrt{N}\mskip\thinmuskip\right\rceil\varphi_i, 2\pi\right).
    \end{aligned}
\end{equation}
We consider two shapes for the endothelium. The first, $y=1\um$, emulates the flat wall typically used in
near-wall simulations of RBCs or platelets. The other attempts to recreate the elongated endothelial cell shape
typical of exposure to high-shear conditions,
\begin{equation}
    y(\theta, \varphi) = 0.75\um + 1\um\cdot\cos^2(\theta-\varphi)\sin^2(\varphi/2).
\end{equation}
The bumps have a prominence of $1\um$. The endothelial surface is raised by $0.75\um$ to avoid it interacting with
the domain boundary. The positions of the surface are chosen to maintain a fixed hematocrit of approximately 34\%
for both endothelial shapes.

As a preliminary validation of the platelet model and to establish baseline platelet motion, we consider two
platelets along a flat wall. They are placed parallel to the wall at distances of $0.3\um$ and $0.5\um$. The
domain does not contain any RBCs. At a distance of $0.3\um$, the platelet is expected to ``wobble\qend,'' a motion in which
the platelet tilts slightly upward and downward, periodically~\cite{King:2005fv}. On the other hand, the platelet
initially $0.5\um$ from the wall should tumble end-over-end. We observe wobbling at a frequency of approximately
$10\persec$ and tumbling with a frequency of approximately $30\persec$. These figures are in reasonable agreement
with other studies~\cite{King:2005fv}. We also note that the edge of the tumbling platelet remains pointed towards
the wall for only 3--$4\ms$.

\subsubsection{Initialization}\label{sec:blood-init}

We assume that the platelets have already been marginated by the RBCs. We think of the domain as having three
layers with the endothelium at the bottom, RBCs on top, and platelets in between. We begin by settling the
endothelium and RBCs before placing platelets in the space between the RBCs and the endothelium.

In addition to the endothelium, we place 2 rows of 4 RBCs, each in their reference configuration, in the domain
with the RBCs' centers of mass on the plane $y=6\um$.  Because the domain is not wide enough to accommodate two
reference RBCs alongside one another, the cells are staggered by $2\um$. These locations are then randomly
translated and rotated while maintaining a distance of at least $2h$ between cells.

Before placing any platelets in the domain, we allow the flow to develop with only the endothelium and RBCs. We
allow the initialization to continue until at least $17\ms$, which is approximately when the first RBC overtakes
its neighbor. From here, we choose a series of times, sampled from a Poisson distribution to be approximately
$3\ms$ apart, at which to begin simulations with platelets.

The RBCs for each of the chosen starting configurations are considerably and unpredictably deformed and have left
a space of a few microns above the endothelium in which we place platelets. To find reasonable starting
orientations for the platelets, we randomly choose points on the endothelium and one on each platelet surface.
Each of the simulations in the upcoming sections contains two platelets, so we choose two points on the
endothelium that are at least $3.9\um$, a platelet diameter plus $4h$, apart. The resulting platelets are spaced
far enough apart to not intersect. We compute normal vectors on the surfaces of the endothelium and platelets at
these points. We orient the platelets so that the normal emanating from the platelet opposes the normal at the
corresponding point on the endothelium. The platelet is then placed so its chosen surface point is separated from
the endothelium point by a random distance between $0.3\um$ and $1\um$. If the generated orientation does not pass
within $0.4\um$ of an RBC, the platelet is accepted. Otherwise, we try again with a different platelet point. This
algorithm typically succeeds within 2 attempts.

This initialization process is performed once for each endothelial configuration and we take the first four
acceptable initial configurations for each. In the following section, we present behaviors found in these
simulations. As a point of comparison, we also consider the same initial configurations for the bumpy wall with
the RBCs removed from the domain.

\subsubsection{Characterization of flow and cell behaviors}

In this section, we catalog the differences in the flow between whole blood along a bumpy and flat wall, and
between flow along a bumpy wall with and without RBCs. We aim to compare the interactions platelets have with RBCs
and the endothelium for these test cases.

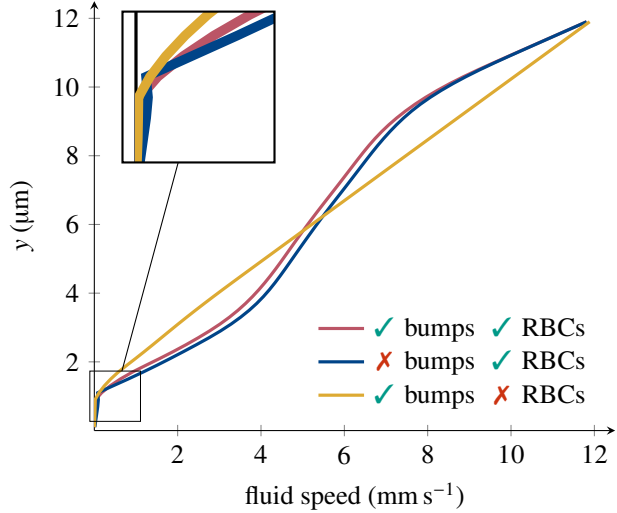
\begin{figure}[tb]
    \centering
    \begin{tikzpicture}[spy using outlines={rectangle, magnification=3,connect spies}]
        \begin{axis}[
                axis lines=center,
                xmin=0,
                xmax=12.5,
                ymin=0,
                ymax=12.5,
                ylabel={$y$ ($\um$)},
                xlabel={fluid speed ($\mmpersec$)},
                xlabel near ticks,
                ylabel near ticks,
                legend pos=south east,
                legend style={draw=none}
            ]
            \addplot[color=tol/contrast/red, very thick, x filter/.code={\pgfmathparse{\pgfmathresult*10}\pgfmathresult}] table [x index=1, y index=0] {rpefast1.prof.dat};
            \addlegendentry{\cmark~bumps\hspace{0.5em}\cmark~RBCs};
            \addplot[color=tol/contrast/blue, very thick, x filter/.code={\pgfmathparse{\pgfmathresult*10}\pgfmathresult}] table [x index=1, y index=0] {rpeflat1.prof.dat};
            \addlegendentry{\xmark~bumps\hspace{0.5em}\cmark~RBCs};
            \addplot[color=tol/contrast/yellow, very thick, x filter/.code={\pgfmathparse{\pgfmathresult*10}\pgfmathresult}] table [x index=1, y index=0] {pefast1.prof.dat};
            \addlegendentry{\cmark~bumps\hspace{0.5em}\xmark~RBCs};

            \coordinate (spypoint) at (axis cs: 0.5, 1);
            \coordinate (spyviewer) at (axis cs: 2.5, 10);
            \spy[width=2cm,height=2cm] on (spypoint) in node [fill=white] at (spyviewer);
        \end{axis}
    \end{tikzpicture}
    \caption[A comparison of fluid velocity profiles]{%
Time- and space-averaged fluid velocity profiles for each of the test cases. The
inclusion of RBCs (red and blue curves) causes the region inhabited by platelets,
1--$4\um$, to experience a higher shear rate than it would without RBCs (yellow curve).
    }\label{fig:flow-profiles}
\end{figure}

Flow profiles are shown in \cref{fig:flow-profiles}. The most notable difference among the three flow profiles is
the nearly Couette flow when RBCs are absent. The only distinction between this and Couette flow with
$\shear=1000\persec$ is the smoother transition at the wall due to the bumps. This is also the distinguishing
feature between the profiles corresponding to bumpy and flat walls in the presence of RBCs. The smooth transition
from the bumps results in marginally slower flow speeds throughout the domain, compared to the flat wall. The
inclusion of RBCs causes the bends in the red and blue curves around $y = 3\um$ and $y = 9\um$. Platelets located
between $y=1\um$ and $y=3\um$ experience higher shear rates in simulations featuring RBCs than in those
without RBCs. The region near the upper boundary also experiences an increased shear rate. In pressure-driven flow
through a tube, we expect parabolic flow. The bend near $y = 9\um$ is therefore nonphysical and arises from
satisfying boundary conditions at the top boundary. However, the increased shear rate in the region between $y =
9\um$ and $12\um$ seems to be useful in deterring RBCs from approaching the upper boundary. This region acts as
another RBC-free layer where the RBCs excurse only infrequently. An exclusionary region of just 1--$2\um$ along
the top boundary increases the effective hematocrit to 37--41\%. Furthermore, the reduced shear rate in the region
containing RBCs results in slower tank-treading compared to a dilute suspension of RBCs, with one period now
lasting approximately $40\ms$.

RBCs are effective at preventing the platelets from moving too far from the endothelium. The furthest observed
distance from the endothelium any platelet takes is just under $1.5\um$. Likewise, RBCs infrequently enter the
RBC-free layer, with some notable exceptions, discussed below. We do not observe any platelet wobbling. Instead,
platelets transiently follow the curve of the bumpy walls, tilt down into the valleys between bumps, and tumble.
Nothing suggests that bumps in the surface of the endothelium alone can sequester platelets, nor do we directly
observe stagnation zones.

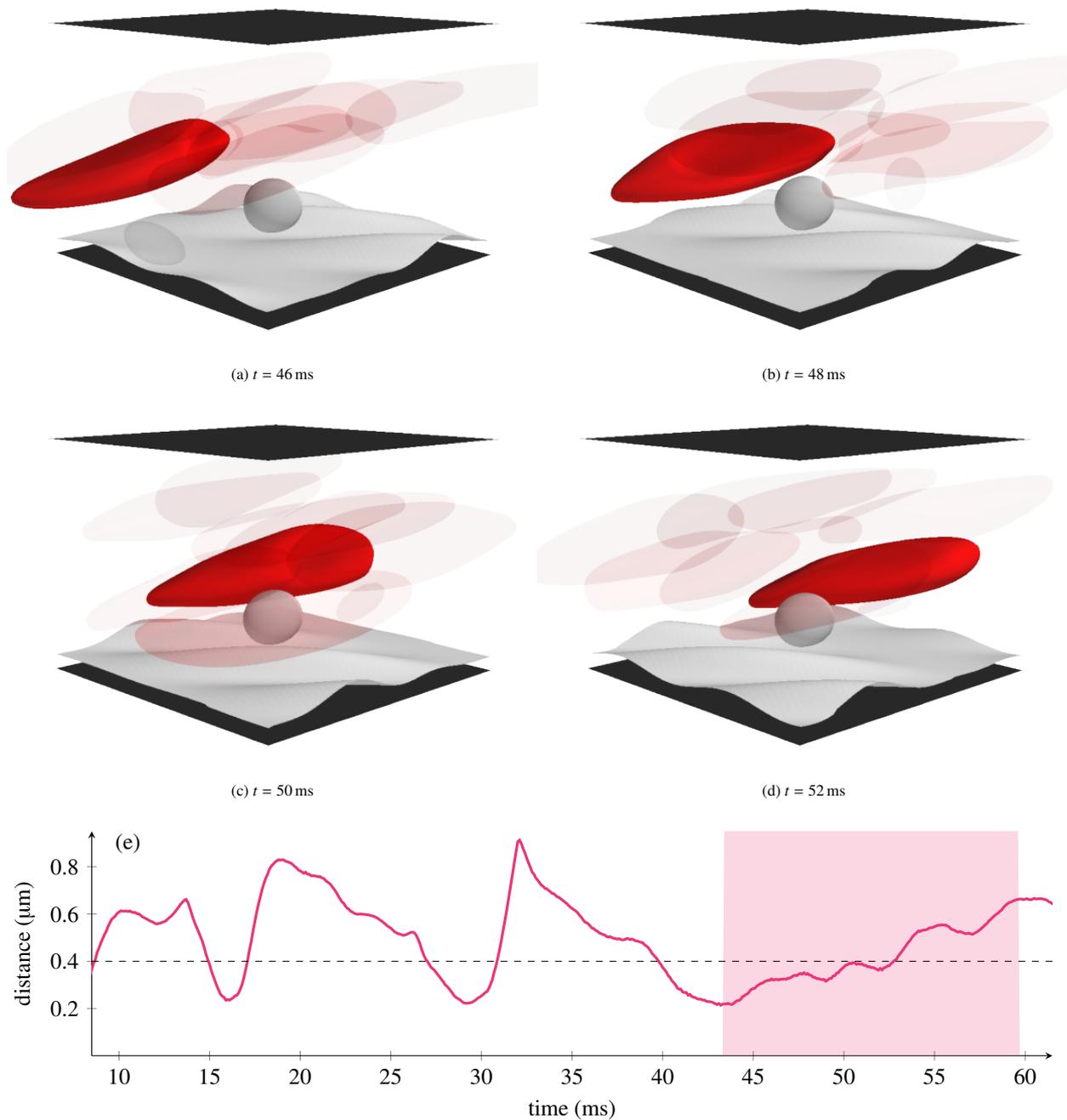
\begin{figure}[tbp]
    \begin{subfigure}[t]{0.5\textwidth}
        \includegraphics[trim=50 75 50 125, clip, width=\textwidth]{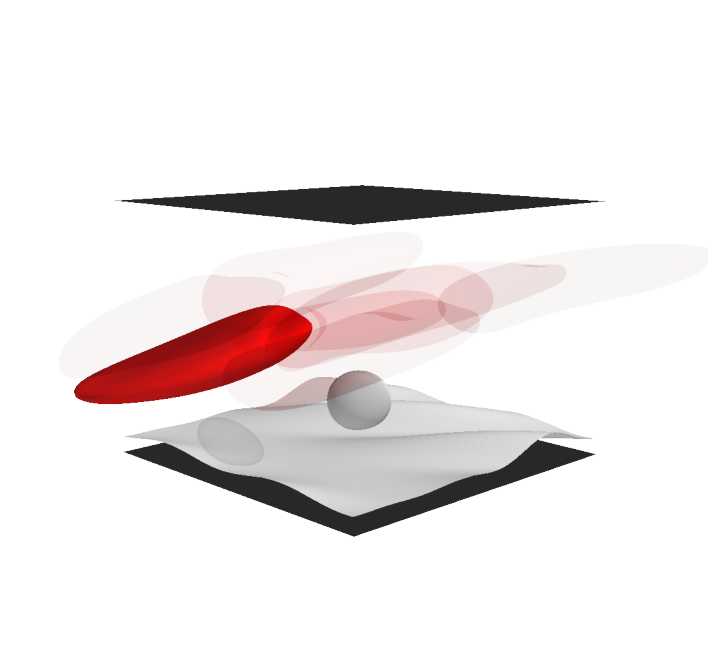}
        \subcaption{$t = 46\ms$}
    \end{subfigure}%
    \begin{subfigure}[t]{0.5\textwidth}
        \includegraphics[trim=50 75 50 125, clip, width=\textwidth]{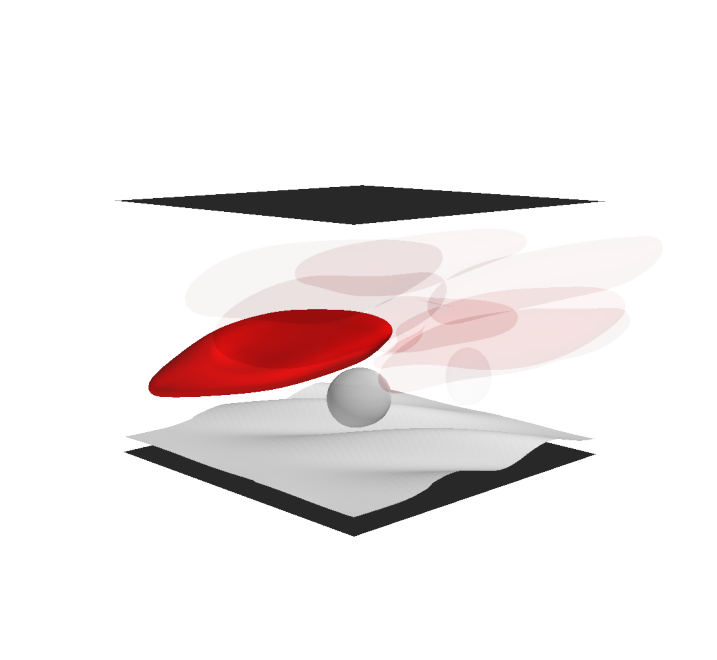}
        \subcaption{$t = 48\ms$}
    \end{subfigure}

    \vspace{11pt}

    \begin{subfigure}[t]{0.5\textwidth}
        \includegraphics[trim=50 75 50 125, clip, width=\textwidth]{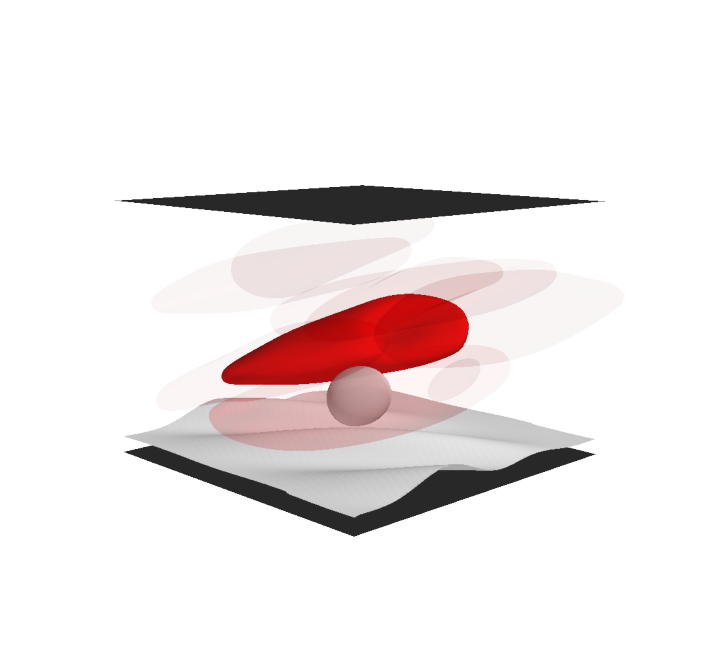}%
        \subcaption{$t = 50\ms$}
    \end{subfigure}%
    \begin{subfigure}[t]{0.5\textwidth}
        \includegraphics[trim=50 75 50 125, clip, width=\textwidth]{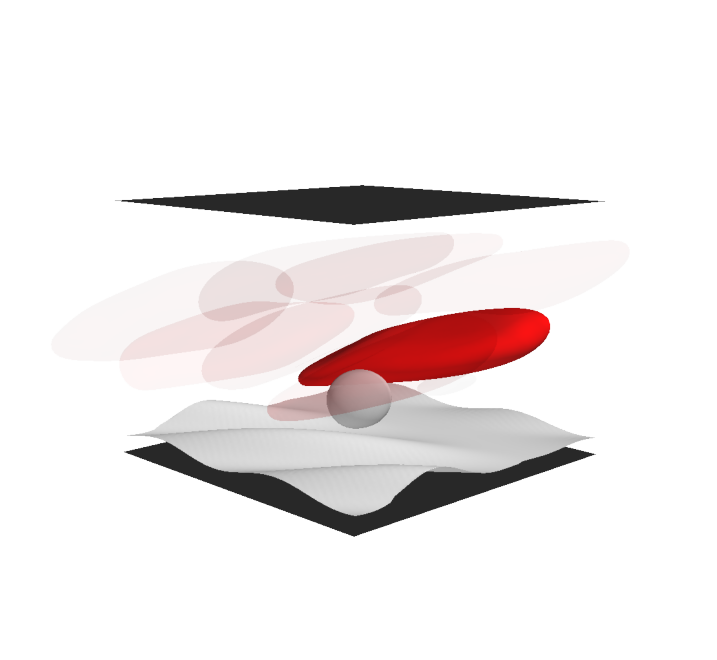}%
        \subcaption{$t = 52\ms$}
    \end{subfigure}

    \vspace{11pt}

    \begin{subfigure}[t]{\textwidth}
        \begin{tikzpicture}
            \begin{axis}[
                width=\textwidth,
                height=2in,
                axis lines=center,
                xmin=8.5,
                xmax=61.5,
                ymin=0,
                ymax=0.95,
                ylabel={distance ($\um$)},
                xlabel={time ($\ms$)},
                xlabel near ticks,
                ylabel near ticks
            ]
                \addplot[color=tol/vibrant/magenta, very thick] table [x index=0, y index=6] {rpefast0.dat};
                \addplot[no marks, color=black, dashed] coordinates {(8.5, 0.4)
                                      (61.5, 0.4)};
                \path[name path=axis] (axis cs: 8.5, 0) -- (axis cs: 61.5, 0);
                \addplot[opacity=0, name path=unicycle] table [x index=0, y index=4] {rpefast0.dat};
                \addplot[fill=tol/vibrant/magenta, fill opacity=0.2] fill between[of=unicycle and axis];

                \node at (axis cs: 10.5, 0.9) {(e)};
            \end{axis}
        \end{tikzpicture}
    \end{subfigure}
    \caption[Platelet unicycling behavior]{%
(a)--(d) Snapshots of a platelet rolling on its edge (``unicycling'') with RBCs, one translucent, flanking either
side. The camera tracks the opaque platelet. The platelet's motion is indicated by the endothelium moving from
right to left. (e) The distance between the platelet and the endothelium. The shaded region indicates that the
orientation of the platelet's short axis is within $45^\circ$ of the vorticity direction. The black dashed line
indicates $2h$ and is the maximum distance that might be considered contact with the endothelium. See~%
\ref{sec:supp} for a video corresponding to this simulation.
    }\label{fig:unicycle}
\end{figure}

Bumpy endothelium simulations without RBCs mimic those with a flat wall; platelets move away from the wall to a
point where they are free to tumble. Unsurprisingly, we observe platelet tumbling for both flat and bumpy walls
with RBCs as well. In Stokes-like flow, a rigid platelet would tumble faster in flow with a higher shear rate. We
might therefore expect the platelet to tumble faster with RBCs. However, RBCs can significantly disturb the fluid
around a platelet, speeding up its motion, slowing it down, or preventing a tumble altogether.

We observe platelets rolling in the flow direction along their edge. Because the platelet in this arrangement is
aligned vertically, part of the edge stays in near-contact with the endothelium while the opposite edge extends
into the region occupied by RBCs.  Contact with RBCs is frequent. These contacts can have a destabilizing effect,
but may also prolong the rolling. \Cref{fig:unicycle}(a)--(d) consists of a series of snapshots illustrating the
behavior. The platelet in this case is flanked by two RBCs, so it does not have the space to topple over until the
RBCs pass a few milliseconds later.

While \cref{fig:unicycle} shows this phenomenon on a bumpy wall, it can also occur above a flat wall. We first
observed this motion with a flat wall, and there it lasted over $40\ms$. The motion was maintained, in part, by an
RBC that rode along the top of the platelet, partially enveloping that platelet. We say that a platelet rolling on
its edge is \emph{unicycling}. \Cref{fig:unicycle}(e) illustrates that the platelet spends more time in contact or
near-contact with the endothelium while unicycling compared to the tumbles near $t=16\ms$ and $t=29\ms$. Though
RBCs seem to control the duration of the unicycling, they are not strictly necessary for unicycling to occur. In
tests with a bumpy wall without RBCs, unicycling is initiated when a platelet rolls sideways, relative to the flow
direction, off of a bump. Without RBCs, the platelet maintains the vertical alignment for a majority of the
simulation thereafter. However, without frequent interaction with RBCs, the platelet in these simulations move
away from the wall. Moreover, while we have not observed it directly, we expect that a lone platelet traveling
over a flat wall would also exhibit unicycling, given the right initial orientation.

\begin{figure}[tbp]
    \begin{subfigure}[t]{0.5\textwidth}
        \includegraphics[trim=50 75 50 125, clip, width=\textwidth]{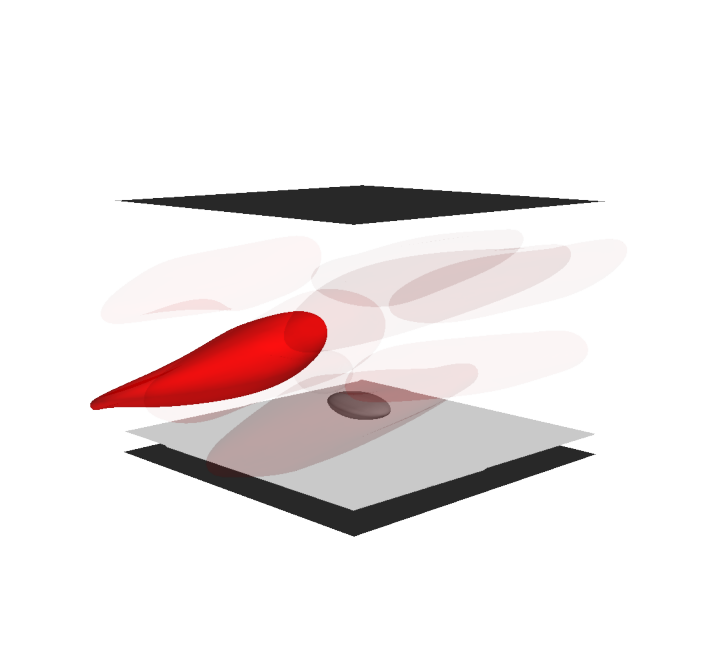}%
        \subcaption{$t = 46.0\ms$}
    \end{subfigure}%
    \begin{subfigure}[t]{0.5\textwidth}
        \includegraphics[trim=50 75 50 125, clip, width=\textwidth]{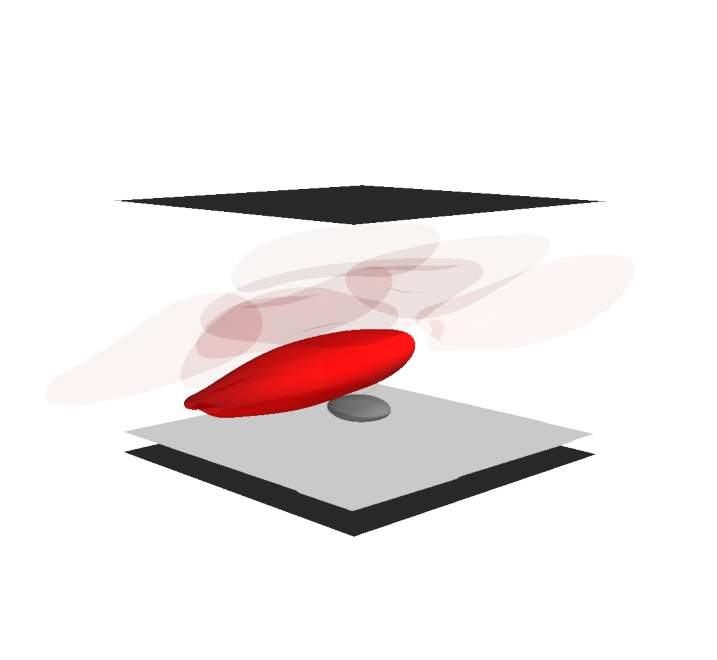}
        \subcaption{$t = 49.5\ms$}
    \end{subfigure}

    \vspace{11pt}

    \begin{subfigure}[t]{0.5\textwidth}
        \includegraphics[trim=50 75 50 125, clip, width=\textwidth]{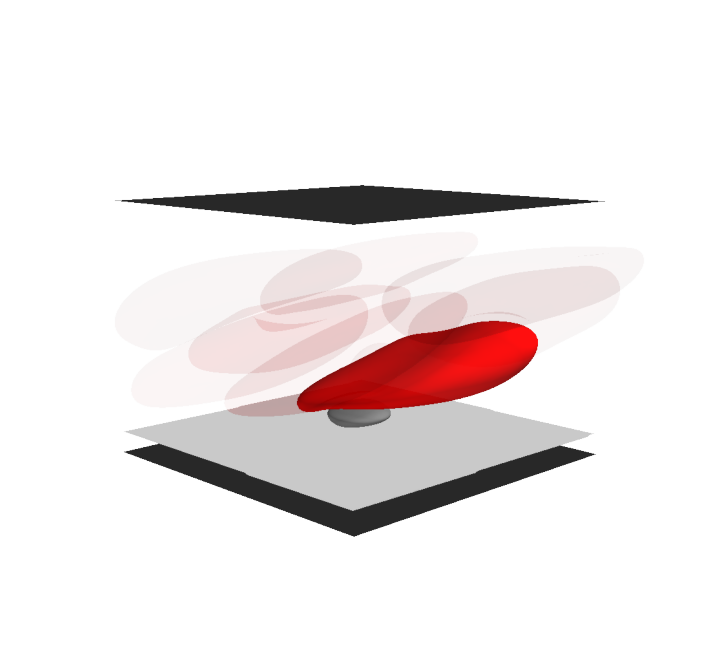}%
        \subcaption{$t = 53.0\ms$}
    \end{subfigure}%
    \begin{subfigure}[t]{0.5\textwidth}
        \includegraphics[trim=50 75 50 125, clip, width=\textwidth]{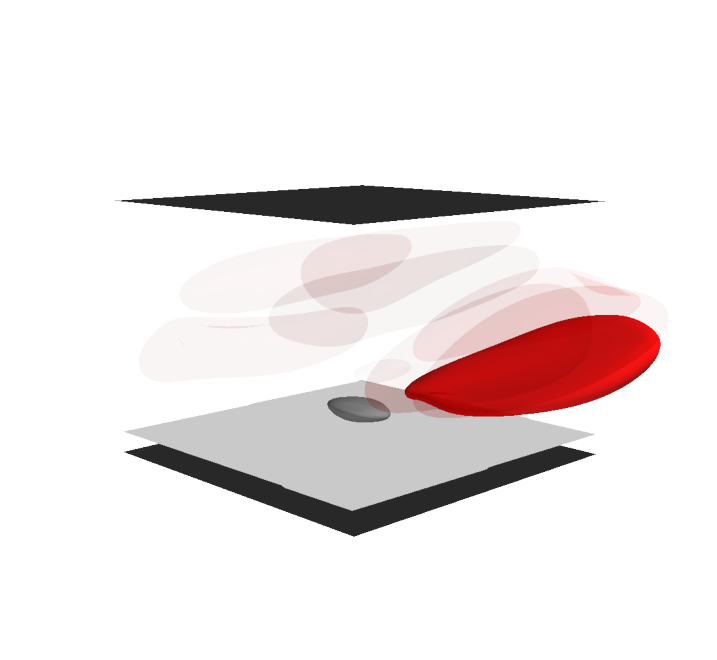}%
        \subcaption{$t = 56.5\ms$}
    \end{subfigure}

    \vspace{11pt}

    \begin{subfigure}[t]{\textwidth}
        \begin{tikzpicture}
            \begin{axis}[
                width=\textwidth,
                height=2in,
                axis lines=center,
                xmin=8.5,
                xmax=61.5,
                ymin=0,
                ymax=3.25,
                ylabel={velocity ($\mmpersec$)},
                xlabel={time ($\ms$)},
                xlabel near ticks,
                ylabel near ticks
            ]
                \addplot[color=tol/vibrant/magenta, very thick] table [x index=0, y index=3] {rpeflat3.vel.dat};
                \path[name path=axis] (axis cs: 8.5, 0) -- (axis cs: 61.5, 0);
                \addplot[opacity=0, name path=deformed, y filter/.code={\pgfmathparse{#1*3}\pgfmathresult}] table [x index=0, y index=4] {rpeflat3.vel.dat};
                \addplot[fill=tol/vibrant/magenta, fill opacity=0.2] fill between[of=deformed and axis];

                \node at (axis cs: 10.5, 3) {(e)};
            \end{axis}
        \end{tikzpicture}
    \end{subfigure}
    \caption[RBC-mediated platelet-endothelial collision]{%
(a)--(d) Snapshots of RBC-mediated collision between a platelet and the endothelium. (a) The platelet attempts to
tumble. (b) and (c) An RBC comes into proximity with the platelet, deflects to avoid the platelet, and pushes the
platelet into the endothelium, thereby preventing the platelet from tumbling. (d) The platelet is free to tumble
again.  (e) The minimum velocity on the surface of the platelet. The shaded region indicates that the relative
change in aspect ratio of the major axes exceeds 4\%. See~\ref{sec:supp} for a video corresponding to this
simulation.
    }\label{fig:rbc-plt-endo-collision}
\end{figure}
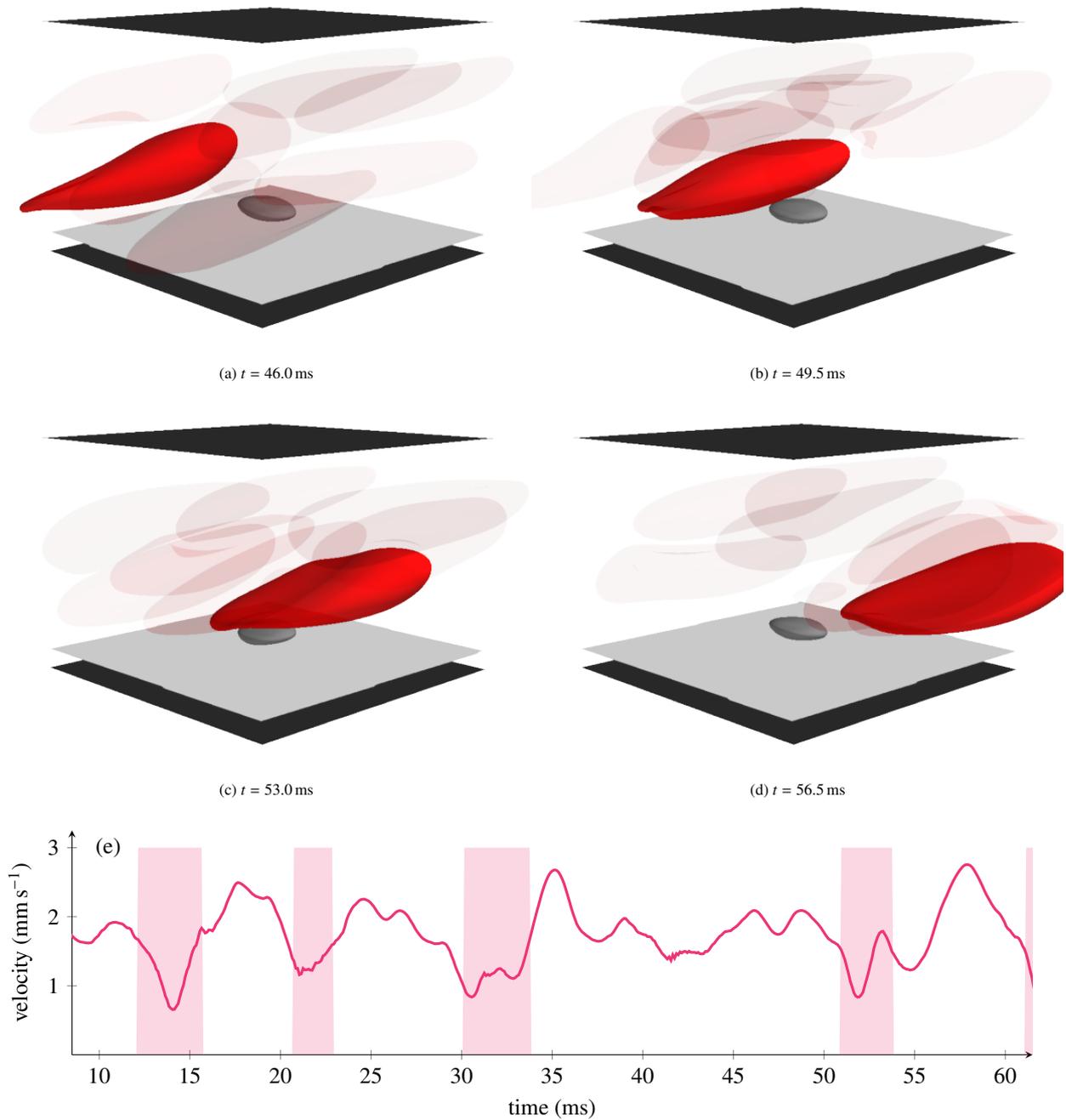

We notice that in simulations with a bumpy endothelium, platelets collide with the bumps. This tends to occur
while the platelet is tumbling, and the edge of the platelet makes contact with the endothelium. This interaction
is characterized by deformations that flatten the edge of the platelet and a \midtilde5\% relative change in the
aspect ratio of the major axes of the platelet. However, the collision need not occur along the edge of the
platelet, nor, indeed, against a bump in the endothelium. Somewhat surprisingly, collisions with a flat wall occur
at roughly the same frequency, suggesting that RBCs mediate this behavior. A clear case of this is illustrated in
\cref{fig:rbc-plt-endo-collision}(a)--(d). We also note that the few milliseconds preceding the unicycling in
\cref{fig:unicycle} correspond to a collision with the wall, showing that this is yet another trigger for
unicycling to occur.

Because the platelet comes into contact with the endothelium, or nearly does so, the platelet slows along the area
of contact. \cref{fig:rbc-plt-endo-collision}(e) shows the correlation between the relative change in aspect ratio
and the reduction in minimum platelet surface velocity. Though the aspect ratio of the platelet changes somewhat
while normally tumbling, changes of 5\% or greater seem to always correspond to interactions with the endothelium.
Collisions with the RBC, for example, result primarily in deformation of the RBC and deflection of the platelet,
which is otherwise relatively unperturbed.

\section{Discussion}\label{sec:conclusion}

In this article, we developed a coherent numerical framework based on the RBF-IB method for whole blood
simulation. We have shown that a continuous energy RBF-based model RBC incorporating dissipative forces exhibits
the traditional tumbling and tank-treading behaviors. We simulated the flow of whole blood involving RBCs and
platelets over a model endothelium. We considered both flat and bumpy endothelial shapes and flow with and without
RBCs along a bumpy wall.

The most prominent result of the whole blood simulations is that simulations involving platelets but neglecting
the influence of RBCs cannot capture the true nature of platelet motion. In fact, platelets near the endothelium
experience augmented shear rates due to RBCs and RBCs act to confine platelets to the cell-free layer. Simulations
without RBCs fail to capture the more irregular aspects of platelet motion. In fact, our simulations show that
interaction with RBCs can disturb otherwise regular wobbling motions exhibited by platelets, and also appear to
delay tumbling. Our results also demonstrate that certain platelet behaviors can only be observed by considering
numerous starting configurations in developed flows involving RBCs.

Another prominent finding is that the effect of RBCs generally overwhelms the effects of the wall topography.
\Cref{fig:flow-profiles} shows that the flat wall yields slightly faster fluid velocities, but the flow profiles
for flat and bumpy walls are qualitatively the same, except for the immediate vicinity of the wall. The result is
a region of space between the bumps with a velocity gradient. Platelets following the shape of the bumpy wall
crest the bump, dip into the region of lower velocity, and tumble. This feature is absent from the flat wall, but
tumbling is not an extraordinary behavior. In either case, we observe ``unicycling\qend,'' a unique behavior in
which the platelet rolls in the flow direction along its edge. Unicycling can be stabilized by RBCs which flank
the platelet, or by an RBC that partially encapsulates the platelet while passing over it. Conversely, it can also
be destabilized by an RBC passing on one side. The endothelial protrusions are also sufficient to orient a platelet
while unicycling, but without RBCs to confine the platelet near the wall, the platelet does not roll along the
wall for long. However, we find that unicycling itself seems to be stable. Unicycling also highlights the need for
3D simulations---it is a behavior that cannot be captured by or predicted from a 2D simulation. We also observe
platelet-endothelial interactions for both endothelial shapes. These interactions are typically caused by RBCs
driving the platelet into the endothelium. The collisions are characterized by significant deformation to the
platelet and a reduction in its speed. From a qualitative standpoint, the endothelial shape alone has minimal
impact on the motion of the platelets, meaning that for modeling flow of whole blood above a healthy blood vessel, a flat wall suffices.

It is also reasonable to consider an alternative interpretation of the flat wall: a model exposed subendothelium.
Near contact with the subendothelium increases the chance of platelet activation.  Unicycling keeps an edge of the
platelet near the wall, without hindering mobility. We see that the vertical alignment is often maintained much
longer than wall contacts from tumbling, and we do not observe wobbling in the presence of RBCs.  We propose
unicycling as an effective means by which platelets survey the vasculature for injury. Of course, the platelet can
only distinguish a healthy vessel from an injury by encountering the necessary chemical signals. Until then, the
platelets unicycle around bumps along the healthy endothelium as well. This also implies that unicycling
indirectly assists in platelet activation for these shear rates.

We also consider yet another alternative interpretation of the bumpy wall. Because the bumps are approximately the
same size as a platelet, we can consider this to be a rough model of a subendothelium with a few deposited
platelets. Under this interpretation, we view platelet-wall contact as interactions between an unactivated
platelet and the subendothelium, when contact occurs in a valley, or between an unactivated platelet and an
activated platelet adhering to the subendothelium, when contact occurs on a bump.  These interactions correlate
with a reduction in velocity at the contact zone and the platelet membrane becomes flattened at the point of
contact. We suggest that the decreased velocity may be sufficient to allow bonds to form between the platelets or
between the platelet and subendothelium. By flattening, the platelet exposes more surface area at the point of
contact, so that the activation signals are more likely to reach the platelet. The observed velocity reduction
seems insufficient for this, but the resolution of our simulations is also unlikely to allow cells to pass within
bonding distance of one another. Overcoming this limitation we leave as a future direction.

\appendix
\section{Boundary error correction for staggered grids}\label{sec:boundary-correction}

The marker-and-cell (MAC) grid~\cite{Welch:1965jv} is a popular method for fluid simulations. Components of
vector-valued quantities are discretized at the center of the corresponding cell faces and scalar-valued
quantities at the cell center. This staggering is a distinguishing feature of the MAC grid. Staggering avoids the
checkerboard instability that arises from using collocated grids~\cite{Wesseling:2001ci}. However, in domains
with non-periodic boundaries, this means that some vector components will encounter situations where satisfying
boundary conditions with linear ghost value extrapolation leads to numerical error. This appendix explores these
errors and provides a resolution that maintains compatibility with the conjugate gradients method.

For a rectangular domain with Dirichlet boundary conditions along its top and bottom and periodic boundaries
elsewhere, consider a region of the domain adjacent to a boundary.  The horizontal component of the fluid
velocity, $u$, is discretized at locations staggered $h/2$ vertically above the bottom of a grid cell, or below
the top of a grid cell, as seen in \cref{fig:discretization}(b). We have boundary data for $u$ along the boundary
at the same $x$ and $z$ coordinates as the grid points. Using a linear interpolant to fill a value at the ghost
cell $h/2$ below the bottom boundary yields
\begin{equation}\label{eq:ghost}
    u_{i, -1, k} = 2\gamma_{i, k} - u_{i, 0, k},
\end{equation}
where $u_{i, -1, k}$ is the value at the ghost cell, $\gamma_{i, k}$ is the boundary datum, and $u_{i, 1, k}$ is
the value at the grid point within the domain. Using the standard 7-point Laplacian then yields
\begin{equation}\label{eq:disc-lap}
    \begin{aligned}
        (\laplacian_h \arr{u})_{i, 0, k}
        &= (D_{xx}\arr{u})_{i, 0, k} + h^{-2}\left[u_{i, -1, k} - 2u_{i, 0, k} + u_{i, 1, k}\right] + (D_{zz}\arr{u})_{i, 0, k} \\
        &= (D_{xx}\arr{u})_{i, 0, k} + h^{-2}\left[2\gamma_{i, k} - 3u_{i, 0, k} + u_{i, 1, k}\right] + (D_{zz}\arr{u})_{i, 0, k},
    \end{aligned}
\end{equation}
where $\arr{u}$ is the vector of $u_{i, j, k}$ values and $D_{xx}$ and $D_{zz}$ are the standard 3-point discrete
second derivative operators with respect to $x$ and $z$, respectively. Because $D_{xx}\arr{u}$ and $D_{zz}\arr{u}$
only involve values at points within the domain, they are known to approximate their continuous counterparts to
second order. We therefore focus on the remaining term. Replacing $u_{i, j, k}$ with $u(hi, h(j+0.5), h(k+0.5))$
and Taylor expanding about $(hi, 0.5h, h(k+0.5))$ yields
\begin{equation}\label{eq:expansion}
    \begin{aligned}
    h^{-2}\left[2\gamma_{i, k} - 3u_{i, 0, k} + u_{i, 1, k}\right]
    &= h^{-2}\left[2\left(u - \frac{h}{2}u_y + \frac{h^2}{8}u_{yy} - \frac{h^3}{48}u_{yyy}\right) - 3u\right. \\
    &\hphantom{=h^{-2}[}\left.+ \left(u + hu_y + \frac{h^2}{2}u_{yy} - \frac{h^3}{6}u_{yyy}\right)+\mathcal{O}(h^4)\right] \\
    &= \frac34u_{yy} - \frac{5h}{24}u_{yyy} + \mathcal{O}(h^2),
    \end{aligned}
\end{equation}
\latin{i.e.}, the leading coefficient is $3/4$ where we expect 1, meaning we obtain a 0\textsuperscript{th}-order
approximation to $u_{yy}$ near the boundary. The case is the same for the upper boundary.

Correcting this while maintaining symmetry of the $\laplacian_h$ operator is not as easy as using a second order
interpolant for the ghost cell value or scaling the stencil for $D_{yy}$ near the boundary by $4/3$. Either of
these options changes the weight of $u_{i, 1, k}$, which destroys symmetry. Instead, we scale any equation
involving a ghost cell by $3/4$, excluding the bracketed terms in~\eqref{eq:disc-lap}. We define the modified
identity operator $\tilde{I}$, which has a value of $3/4$ on the diagonal for rows corresponding to near-boundary
equations and 1 elsewhere on the diagonal. We define the modified discrete Laplacian
\begin{equation}\label{eq:mod-disc-lap}
    \tilde{\laplacian}_h = \tilde{I}D_{xx} + D_{yy} + \tilde{I}D_{zz},
\end{equation}
where, away from the boundary, $D_{yy}$ is the standard 3-point discrete second derivative with respect to $y$,
and defined according to~\eqref{eq:disc-lap} otherwise.  The resulting $D_{yy}$ is symmetric. In addition to
replacing the discrete Laplacian $\laplacian_h$ with $\tilde{\laplacian}_h$, a timestepping scheme using this
correction must replace instances of the identity matrix, many of which are not written explicity, with the
modified identity $\tilde{I}$. A stage of such a timestepping scheme takes the form
\begin{equation*}
    (\tilde{I} - \alpha k\tilde{\laplacian}_h)\arr{u}^n = (\tilde{I} + (1-\alpha)k\tilde{\laplacian}_h)\arr{u}^{n-1} + B_h\arr{\gamma} + \tilde{I}\arr{f},
\end{equation*}
where $\alpha$ is either 0.5 or 1, depending on the scheme and stage, and $B_h$ is an appropriate operator that
modifies equations near the boundary with boundary data. This corrected scheme involves modified operators, but
requires only an additional diagonal multiplication compared to the uncorrected scheme. The Helmholtz operator on
the left-hand side remains symmetric positive definite after correction, and is therefore suitable for solution
via conjugate gradients.

\section{Force models}\label{sec:forces}

In this section, we list the force densities associated with each of the constitutive laws given in
\cref{sec:energy}. For simplicity, we adopt the Einstein summation notation. A Greek letter featuring as a
subscript and superscript within a term, \latin{e.g.}, $a_\alpha b^\alpha$, indicates summation over $\{1, 2\}$
for that letter. We also adopt the comma notation for partial differentiation, where subscripts following a comma
indicate partial differentiation with respect to the corresponding coordinates, \latin{e.g.},
$\phi_{,\alpha\beta} = \partial^2\phi/\partial\q[\alpha]\partial\q[\beta]$. The surface coordinates $\q[1]$ and
$\q[2]$ correspond to $\theta$ and $\varphi$ of \cref{sec:energy,sec:rbfs} in either order.

We begin with Hookean and damped spring forces, whose linearized force densities together take the form
\begin{equation}
    \F_\text{spring} = -k (\X - \X') - \eta(\U - \U').
\end{equation}
It is common practice to treat each spring individually, so that the quadrature weight $\weight[j]$ is absorbed
into the coefficients: $k$ has units of force per length and $\eta$ units of force-time per length. Implementation
of these forces requires no geometric information outside of positions and velocities. These are only used for the
endothelium.

Next, we consider tension, which generates forces based on stretching or compressing of the elastic surface. This
is somewhat more complicated than the Hookean spring case.  First, we define the \emph{metric tensor},
$\metric_{\alpha\beta} = \X_{,\alpha}\cdot\X_{,\beta}$  (which encodes local information about distance and area)
and its inverse, $\metric^{\alpha\beta}$ (sometimes called the dual metric). We similarly define the metric tensor
for the reference configuration, $\reference\metric_{\alpha\beta}$, and its dual. This allows us to write the
Green-Lagrange strain tensor as
\begin{equation}
    \epsilon_\alpha^\beta = \frac12\left(\metric_{\alpha\mu}\reference\metric^{\mu\beta}-\Kronecker_\alpha^\beta\right),
\end{equation}
The invariants of this tensor encode information about relative changes in lengths and areas, and therefore can be
used to write a tension force density. The invariants can be expressed as
\begin{align}
    I_1 &= 2\epsilon_\mu^\mu, \\
    I_2 &= 4\epsilon + I_1,
\end{align}
where $\epsilon_\mu^\mu$ and $\epsilon$ are the trace and determinant, respectively, of
$\epsilon_\alpha^{\smash\beta}$. For Skalak's Law~\eqref{eq:skalak-law} and neo-Hookean tension~%
\eqref{eq:neohookean}, we first define the second Piola-Kirchhoff stress tensor using the invariants $I_1$ and
$I_2$ as
\begin{equation}
    s^{\alpha\beta} = 2\frac{\partial W}{\partial I_1} \hat{g}^{\alpha\beta} + 2I_2\frac{\partial W}{\partial I_2} g^{\alpha\beta},
\end{equation}
where $W(I_1, I_2)$ is the tension energy density function (see \cref{sec:energy}). This in turn allows us to
define the tension force density~\cite{Maxian:2018ek}
\begin{equation}\label{eq:tension-force}
    \F_\text{tension} = \frac{1}{\sqrt{\reference\metric}}{\left(\sqrt{\reference\metric}s^{\alpha\beta}\X_{,\beta}\right)}_{,\alpha},
\end{equation}
where $\X_{,\beta}$ refer to the tangent vectors on the surface. Because the tension force density is expressed in
relation to the reference configuration, the force is computed by multiplying by quadrature weights for the
reference configuration, which do not change over the course of a simulation.

Our models also contain terms to penalize bending. First, given tangent vectors $\X_{,1}$ and $\X_{,2}$, the unit
normal vector to the surface is given as
\begin{equation}\label{eq:unit-normal}
\n = \frac{1}{\sqrt{\metric}} (\X_{,1}\times\X_{,2}).
\end{equation}
This in turn allows us to define the symmetric tensor $b_{\alpha\beta} = \n\cdot\X_{,\alpha\beta}$, which contains
the coefficients of the second fundamental form, and the \emph{shape tensor}
$K_\alpha^\beta = b_{\alpha\mu}\metric^{\mu\beta}$. The principal curvatures are the eigenvalues of the shape
tensor. More importantly, the trace of this tensor is twice the mean curvature $2H = K_\mu^\mu$, and its
determinant $K$ is the Gaussian curvature. We use $H$ and $K$ within a standard expression for the Canham-Helfrich
force density~\cite{Zhongcan:1989ue} to obtain the following bending force density:
\begin{equation}\label{eq:bending-force}
    \F_\text{CH} = -4\kappa\left(\laplacian(H-H')+2(H-H')(H^2-K+HH')\right)\n,
\end{equation}
where $\laplacian$ is the Laplace-Beltrami operator. We can compute $H$ and $K$ using the formulas above, but
$\laplacian H$ requires up to fourth derivatives of $\X$. In numerical simulations, we compute $H$ pointwise and
apply the discrete Laplace-Beltrami operator to obtain this $\laplacian H$.

Finally, we consider dissipative forces, which cause the membrane to exhibit a viscoelastic response to strain.
With surface velocity $\U$, the metric tensor changes in time according to
\begin{equation}
    \dot{\metric}_{\alpha\beta} = \U_{,\alpha}\cdot\X_{,\beta} + \X_{,\alpha}\cdot\U_{,\beta}.
\end{equation}
The dissipative force density takes the form~\cite{Rangamani:2012hi}
\begin{equation}\label{eq:dissip-force}
    \F_\text{dissip} = \frac{\nu}{\sqrt{\metric}}{\left(\sqrt{\metric}\metric^{\alpha\mu}\dot{\metric}_{\mu\lambda}\metric^{\lambda\beta}\X_{,\beta}\right)}_{,\alpha},
\end{equation}
where $\nu$ is the membrane viscosity.

In general, for all force calculations, it is possible to rewrite~\eqref{eq:tension-force},~%
\eqref{eq:dissip-force}, and the Laplace-Beltrami operator in~\eqref{eq:bending-force} in terms of first and
second derivatives with respect to parametric variables, thereby ensuring that we require only discrete first and
second derivative operators to compute a wide variety of forces.

\section{Supplementary materials}\label{sec:supp}
Videos for Figures~\ref{fig:unicycle} and~\ref{fig:rbc-plt-endo-collision} can be found
in the online supplement.

\section*{Acknowledgments}
The authors thank Dr.\ Robert M.\ Kirby for access to computing resources at SCI Institute. A.\ T.\ K., A.\ B.,
and A.\ L.\ F.\ acknowledge support for this project under NHLBI grant 1U01HL143336 and NSF grants DMS-1716898 and DMS-1521748. V.\ S.\ acknowledges support for this project under NSF grants CISE CCF 1714844 and DMS-1521748.

\bibliography{modeling-paper}
\end{document}